\documentclass[10pt]{article}

\usepackage[twoside]{geometry}
\usepackage{amssymb,amsmath,amsthm,tikz,multirow,nccrules,graphicx,subfig}
\usetikzlibrary{arrows,calc}
\usepackage{indentfirst}
\usepackage{url}
\usepackage[hidelinks, hyperindex]{hyperref}
\usepackage{float}
\usepackage{mathdots}  
\usepackage{enumitem}
\usepackage{booktabs}
\usepackage{longtable}
\usepackage{adjustbox}
\usepackage{dcolumn}
\newcolumntype{L}{>{$}l<{$}}

\title{Non-side-to-side tilings of the sphere by congruent triangles with any irrational angle}
\author{Wen Chen,  Jinjin Liang, Erxiao Wang\thanks{Corresponding author (wang.eric@zjnu.edu.cn).  Research was supported by National Natural Science Foundation of China NSFC-RGC 12361161603 and Key Projects of Zhejiang Natural Science Foundation LZ22A010003.}\\
	Zhejiang Normal University}

\newcommand\bbb{\beta}
\newcommand\ccc{\gamma}
\newcommand{\thin}{\hspace{0.1em}\rule{0.7pt}{0.8em}\hspace{0.1em}}
\newcommand\bn{\boldsymbol{n}}
\newcommand\aaa{\alpha}

\newtheorem{theorem}{Theorem}
\newtheorem{lemma}[theorem]{Lemma}
\newtheorem{remark}[theorem]{Remark}

\newtheorem*{theorem*}{Theorem}
\theoremstyle{definition}
\newtheorem*{definition*}{Definition}
\newtheorem*{case*}{Case}
\numberwithin{equation}{section}

\begin{document}
	\date{}
\maketitle
\begin{abstract}
	We develop the basic and new tools for classifying non-side-to-side tilings of the sphere by congruent triangles. Then we prove that, if the triangle has any irrational angle in degree, such tilings are: a sequence of 1-parameter families of triangles each admitting many 2-layer earth map tilings with $2n$($n\geq3$) tiles, together with rotational  modifications for even $n$; a 1-parameter family of triangles each admitting a unique tiling with $8$ tiles; and a sporadic triangle admitting a unique tiling with $16$ tiles. Then a scheme is outlined to classify the case with all angles being rational in degree, justified by some known and new examples.
	
	
\end{abstract}

\section{Introduction}

A century ago, Sommerville \cite{som} started to study tilings of the sphere by congruent triangles, and Ueno and Agaoka \cite{ua} classified all edge-to-edge cases in 2002 after Davies' 1967 report \cite{dav}. In recent years, a lot of progress has been made on spherical tilings: multihedral tilings by regular polygons of three or more sides are classified in \cite{aehj}; edge-to-edge monohedral tilings by pentagons or quadrilaterals are classified in \cite{wy1, wy2, awy, lwy, cly2, slw, lw1, lw2, lw3, cly1}. However, we still do not know all triangle tilings which are not edge-to-edge, after Dawson and Doyle's many earlier explorations \cite{dawson2001, dawson2003, dawson2006, dawson2006-1, dawson2007}. This paper and a sequel \cite{jwy} aim to complete this classification following Rao's remarkable work \cite{rao}. 

We will use `side-to-side' instead of `edge-to-edge' to emphasize the difference between the corners, sides of a polygonal tile and the vertices, edges of a tiling structure (see \cite{gs,ac1}). Recall that some side-to-side quadrilateral tilings produce new classes of non-side-to-side  triangular ones in our recent work \cite{lw1, lw2, lw3}. Let us list all known examples as five classes below, after introducing some notations.
      \begin{figure}[H]
      	\centering
      	\begin{tikzpicture}
      		\draw
      		(-0.8,-0.8) -- (0,0.8)  -- (0.8,-0.8);

      		\draw[line width=1]
      		(-0.8,-0.8) -- (0.8,-0.8);	    
      		
      		\node at (-0.5,-0.6) {\small $\bbb$};
      		\node at (0.5,-0.6) {\small $\ccc$};
      		\node at (0,0.4) {\small $\aaa$};

      		\node at (0,-1.1) {\small $a$};
      		\node at (0.7,0) {\small $b$};
      		\node at (-0.7,0) {\small $c$};
      		
      	\end{tikzpicture} 
      	\caption{A geodesic triangle on the sphere.}
      	\label{quad}
      \end{figure}
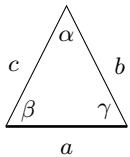

     For a triangle, denote its sides by  $a$, $b$, $c$ and corners by $\aaa$, $\bbb$, $\ccc$ as indicated in Figure \ref{quad}. We use $\aaa^m\bbb^n\ccc^l$ to mean a vertex having $m$ copies of $\aaa$, $n$ copies of $\bbb$ and $l$ copies of $\ccc$. The anglewise vertex combination(s), or AVC, collects all vertices in a tiling. Then the notation $T(6\aaa^4, 6\aaa\ccc^5, 2\bbb^6;\, 6\bbb^3)$ means the tiling has exactly $6$ full vertices $\aaa^4$, $6$ full vertices $\aaa\ccc^5$, $2$ full vertices $\bbb^6$, $6$ half vertices $\bbb^3$, and is uniquely determined by them. Here a full or half vertex means all corners there sum to $2\pi$ or $\pi$ respectively, see Section 2 for more details. A half vertex lies in the interior of a tile's side, and must appear in any non-side-to-side tiling. If the tiling is not uniquely determined by the AVC, then we use  ``$\{12\aaa^{3}\bbb^2,2\bbb^6;\,12\bbb^3\}$: 3" to mean that there exist exactly three different tilings (up to rotations and global flip of the whole sphere) with the same set of vertices. 
       
\subsubsection*{1. Two types of earth map tilings and their  modifications}
\begin{figure}[H]
	\centering
	\includegraphics[width=0.133\textwidth]{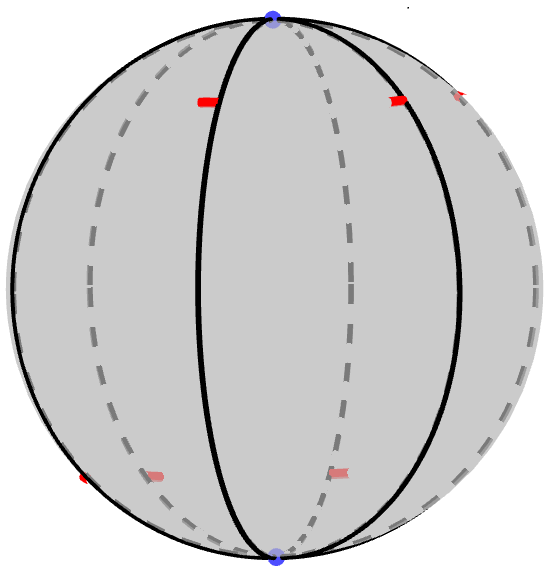}\hspace{0.3cm}
	\includegraphics[width=0.14\textwidth]{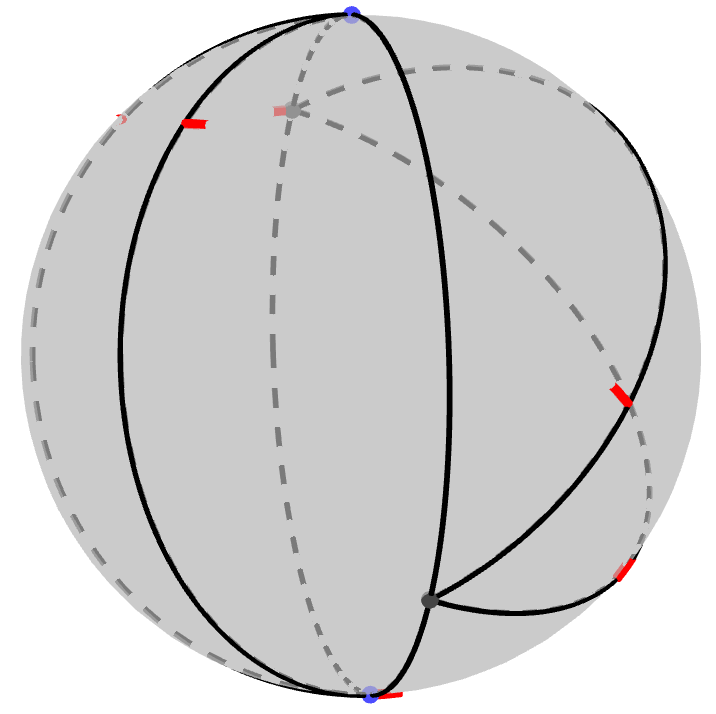}\hspace{0.3cm}
	\includegraphics[width=0.14\textwidth]{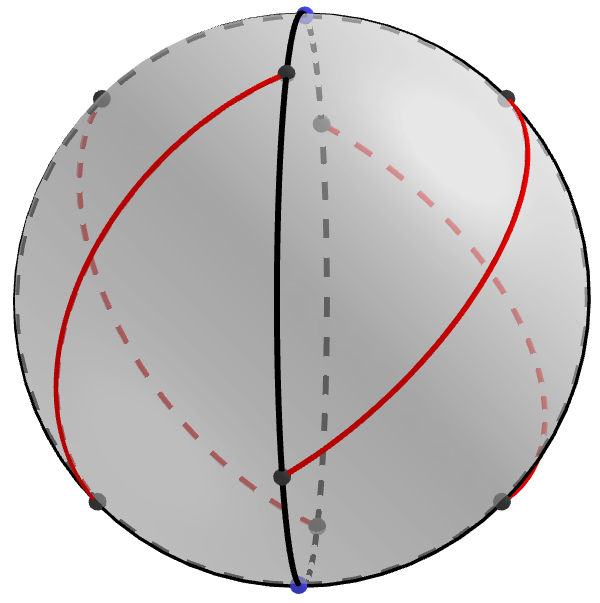}\hspace{0.3cm}
	\includegraphics[width=0.145\textwidth]{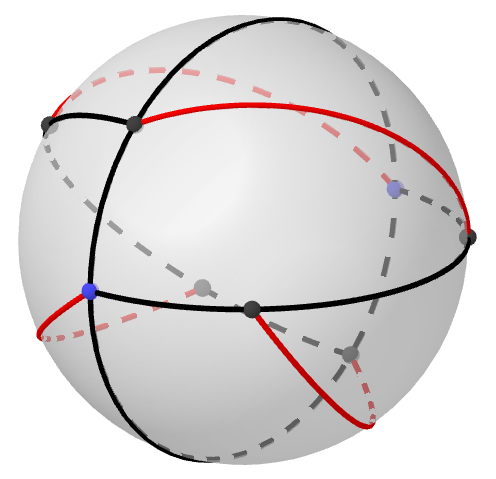}\hspace*{0.3cm}
	\includegraphics[width=0.14\textwidth]{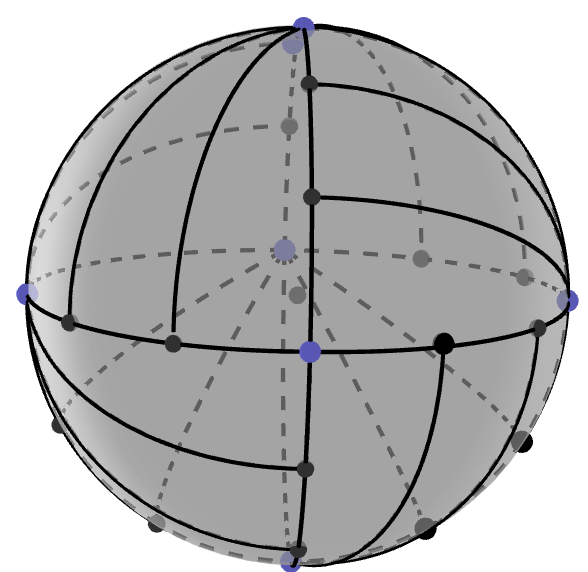}\hspace*{0.3cm}
	\includegraphics[width=0.137\textwidth]{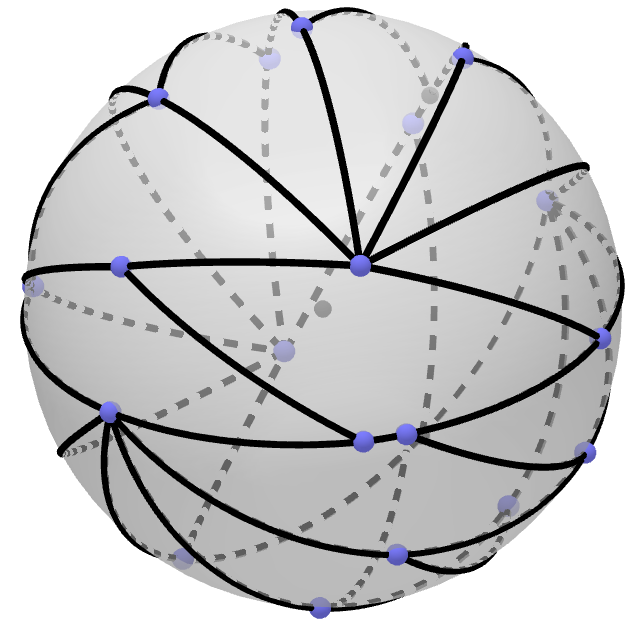}
	\caption{One-layer and two-layer earth map tilings and their  modifications.}
	\label{1.2}
\end{figure}

    Similar to the time zones of the earth, the sphere can be divided into $k\ge2$ congruent $2$-gons or lunes, which may be viewed as degenerate triangles ($\alpha=\pi$). These are called \textbf{one-layer} earth map tilings, and our Lemma \ref{degen} will characterize them and their rotation modifications for even $k$ as all the tilings by congruent degenerate triangles. See the first two pictures of Figure \ref{1.2}. The extended edges are $a=b+c=c+b$ or $a+b+c=2\pi=b+c+b+c=2a$.
    
    Furthermore, each $2$-gon can be divided (with a free parameter) to two congruent triangles to produce \textbf{two-layer} earth map tilings and their modifications, as shown in the last four pictures of Figure \ref{1.2}. The extended edges are $a+b=b+a$, $a+b+a+b=2\pi$, $a=3c=\frac{\pi}{2}$, both $a+b=b+a$ and $a+2b=2b+a$, respectively.

\subsubsection*{2. Various subdivisions of Platonic solids  and their rotation modifications}   
\begin{figure}[H]
     	\centering
     	{\includegraphics[width=0.114\textwidth]{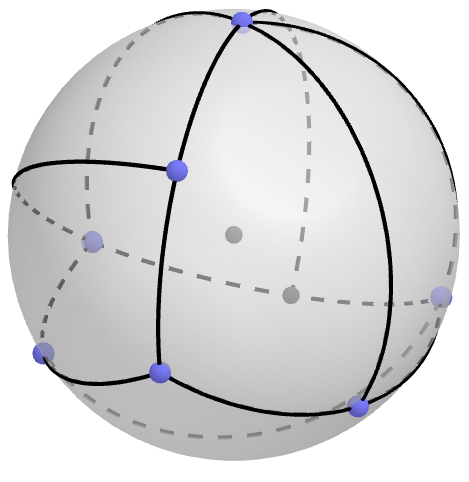}}\hspace*{0.1cm}
\includegraphics[width=0.119\textwidth]{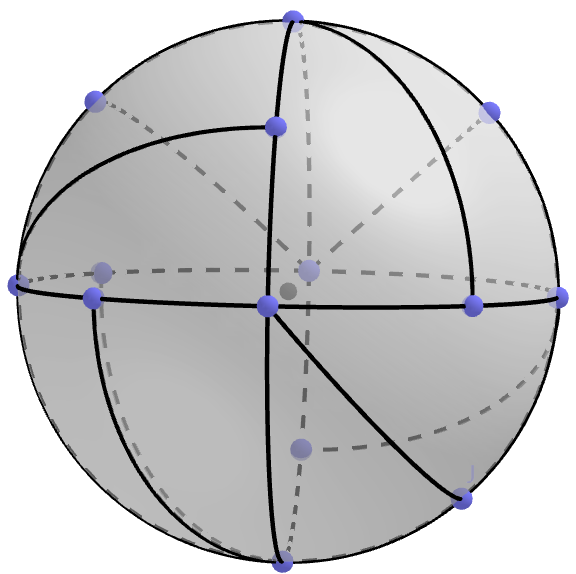}\hspace*{0.1cm}\includegraphics[width=0.117\textwidth]{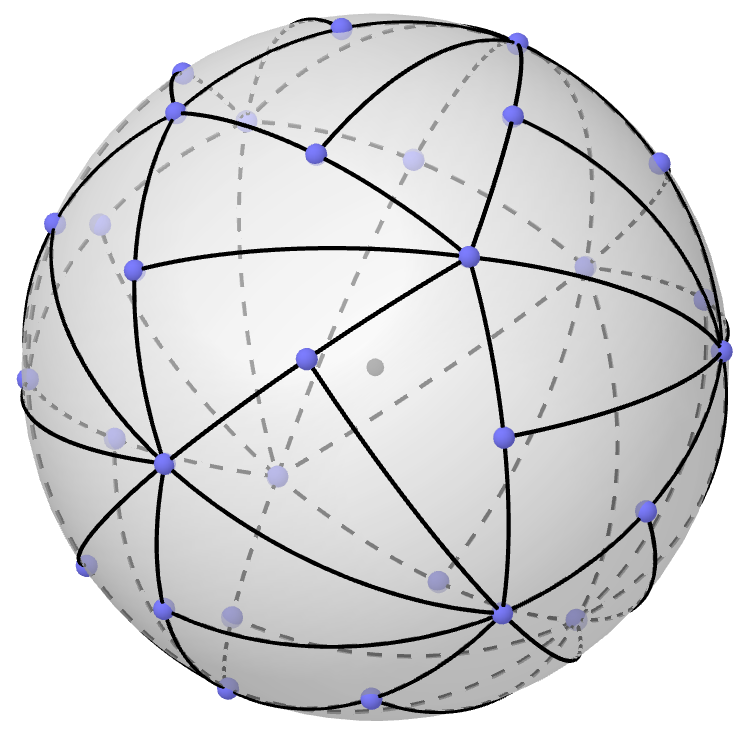}\hspace*{0.1cm}
     	\includegraphics[width=0.117\textwidth]{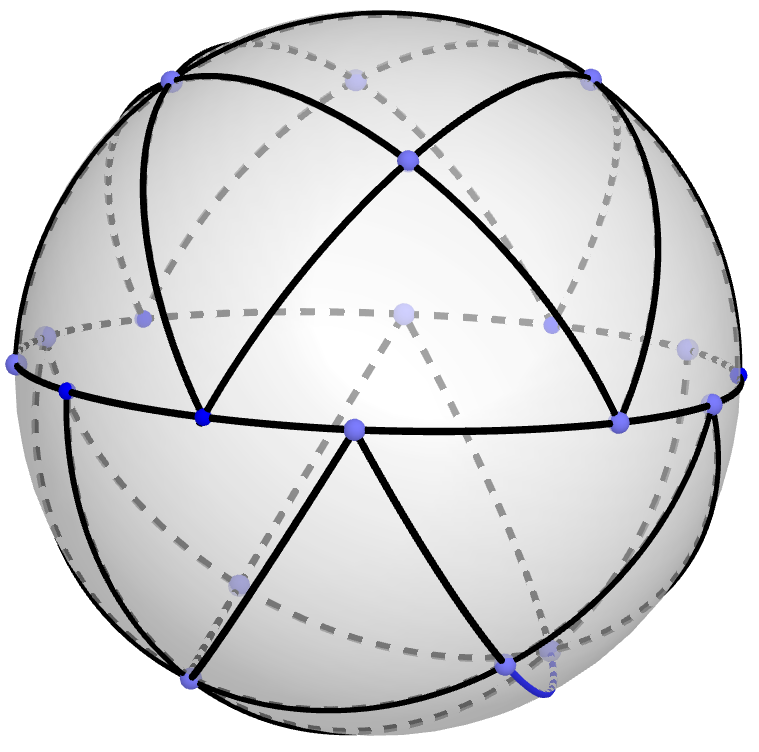}\hspace*{0.1cm}
     	{\includegraphics[width=0.114\textwidth]{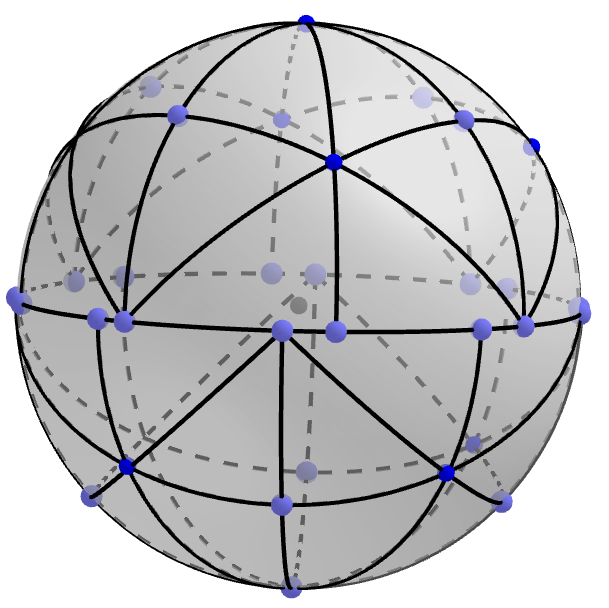}}\hspace*{0.06cm}
     	{\includegraphics[width=0.114\textwidth]{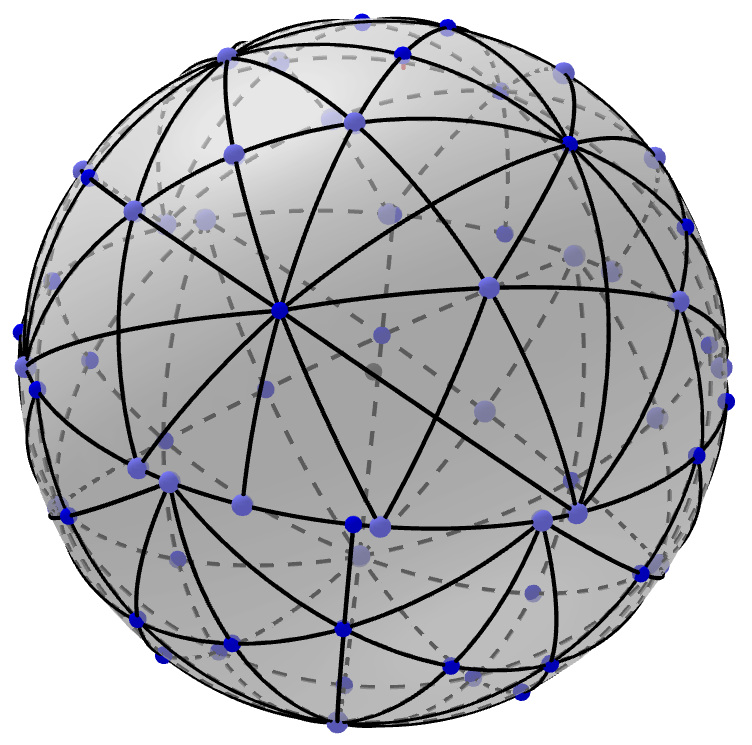}}\hspace*{0.06cm}
     	\includegraphics[width=0.12\textwidth]{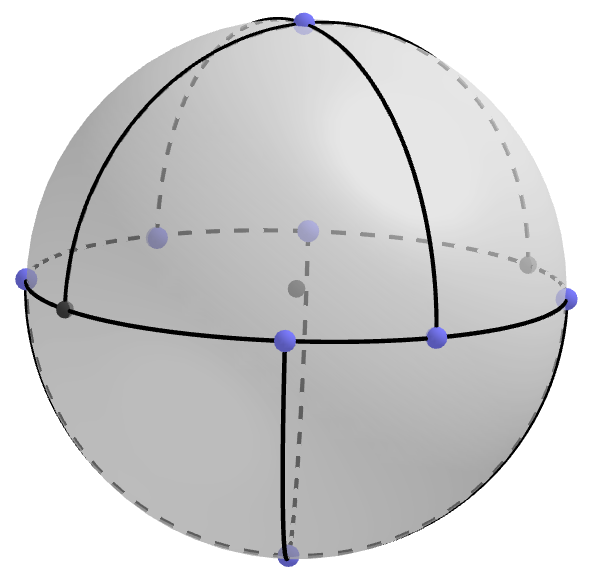}\hspace*{0.06cm}
     	\includegraphics[width=0.114\textwidth]{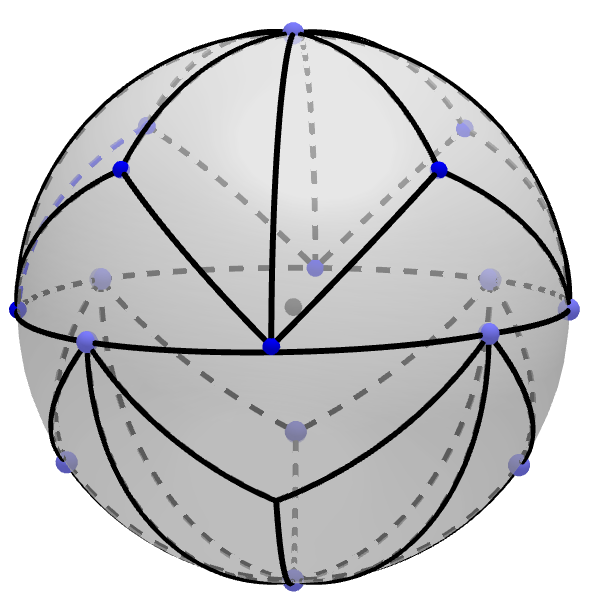}
     	\caption{Triangular subdivisions of Platonic solids and their rotation modifications.}
     	\label{subdivision}
     \end{figure}

Each regular triangular face of the tetrahedron, octahedron and icosahedron can be subdivided into two congruent right triangles in one of the three different ways, to produce many non-side-to-side tilings (see the first three pictures of Figure \ref{subdivision}). The barycentric subdivision of any Platonic solid has many great circles of sides, along which most rotations will produce non-side-to-side tilings (see the next three pictures). The octahedron itself and its central subdivision also have many great circles of sides, along which most rotations will produce non-side-to-side tilings (see the last two pictures).      
     
\subsubsection*{3. Three more types of isosceles triangle monotiles and their tilings} 
\begin{figure}[H]
	\centering
	\includegraphics[width=0.137\textwidth]{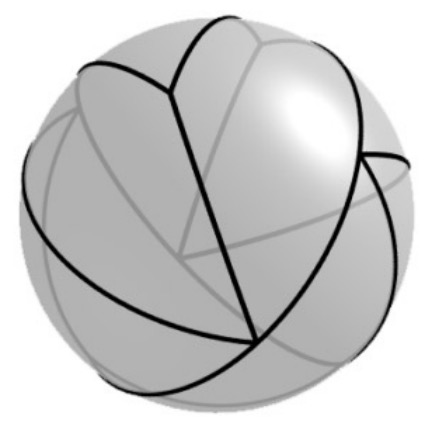}\hspace*{0.5cm}
	\includegraphics[width=0.1344\textwidth]{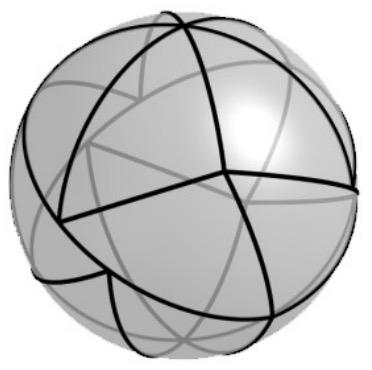}\hspace*{0.5cm}
	\includegraphics[width=0.137\textwidth]{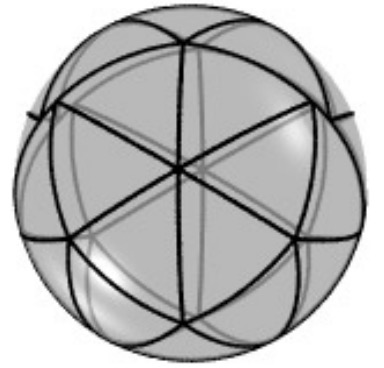}\hspace*{0.5cm}
	\includegraphics[width=0.137\textwidth]{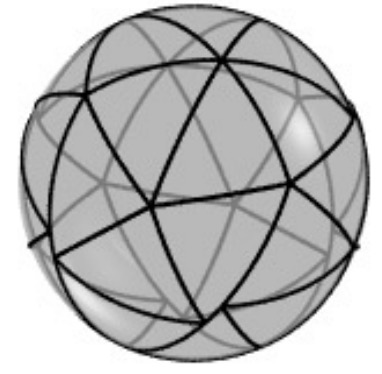}\hspace*{0.5cm}
	\includegraphics[width=0.137\textwidth]{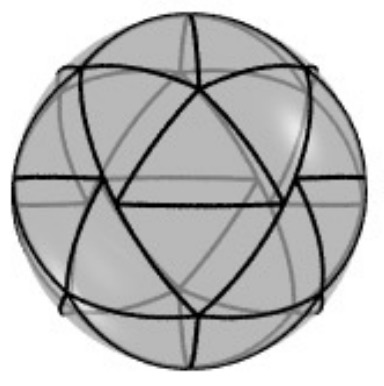}\hspace*{0.5cm}
	\caption{Non-side-to-side tilings by congruent isosceles triangles.}
	\label{iso figure}
\end{figure}

 Non-side-to-side tilings of the sphere by isosceles triangles have been classified by Dawson \cite{dawson2001,dawson2003}. Excluding the triangles mentioned above, there are three additional isosceles triangle monotiles. The 3D pictures for these examples are shown in Figure \ref{iso figure}. The exact geometric data and the extended edges are given in Table \ref{iso}. 

\begin{table}[htp]
	\centering
	\caption{Three more isosceles triangle monotiles and their tilings by Dawson.}
	\label{iso}
	\begin{tabular}{|c|c|c|}
		\hline
		$f$ & isosceles triangle monotiles  &all vertices and extended edges \\
		\hline
		& $(\frac{5}{6}\pi,\frac{1}{3}\pi,\frac{1}{3}\pi)$ &\\
		8&$a=\pi-\arccos(\frac{2\sqrt{3}}{3}-\frac{1}{3})\approx2.5346$& $T(4\aaa^2\bbb;\, 4\bbb^3)$\\
		&$b=\pi-\arccos(\frac{(-2+\sqrt{3})\sqrt{3}}{3})\approx1.4155$&$a+b=b+a$\\
		&$c=\pi-\arccos(\frac{(-2+\sqrt{3})\sqrt{3}}{3})\approx1.4155$&\\
		\hline
		&$(\frac{5}{9}\pi,\frac{1}{3}\pi,\frac{1}{3}\pi)$&\\
		18&$a=\arccos(-\frac{4 \cos(\frac{4 \pi}{9})}{3}+\frac{1}{3})\approx 1.4688$&$T(6\aaa^3\bbb,2\bbb^6;\, 6\bbb^3)$\\
		&$b=\pi-\arccos(\frac{(\cos(\frac{4 \pi}{9})-1) \sqrt{3}}{3 \sin(\frac{4 \pi}{9})})\approx1.0651$& $a+b=b+a$\\
		&$c=\pi-\arccos(\frac{(\cos(\frac{4 \pi}{9})-1) \sqrt{3}}{3 \sin(\frac{4 \pi}{9})})\approx1.0651$&\\
		\hline
		&$(\frac{4}{9}\pi,\frac{1}{3}\pi,\frac{1}{3}\pi)$&\\
		36&$a=\arccos(\frac{4 \cos(\frac{4 \pi}{9})}{3}+\frac{1}{3})\approx0.9705$&$\{12\aaa^{3}\bbb^2,2\bbb^6;\, 12\bbb^3\}$: 3 tilings\\
		&$b=\arccos(\frac{(\cos(\frac{4 \pi}{9})+1) \sqrt{3}}{3 \sin(\frac{4 \pi}{9})})\approx0.8120$&$a+b=b+a$\\
		&$c=\arccos(\frac{(\cos(\frac{4 \pi}{9})+1) \sqrt{3}}{3 \sin(\frac{4 \pi}{9})})\approx0.8120$&\\
		\hline
	\end{tabular}
\end{table}

\subsubsection*{4. Seven more right triangle monotiles and their tilings}
\begin{figure}[H]
	\centering
	\includegraphics[width=0.133\textwidth]{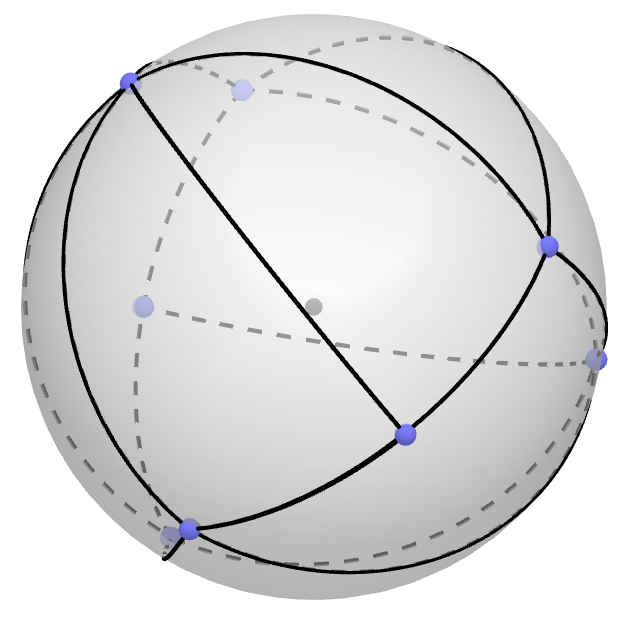}\hspace*{0.1cm}
	\includegraphics[width=0.133\textwidth]{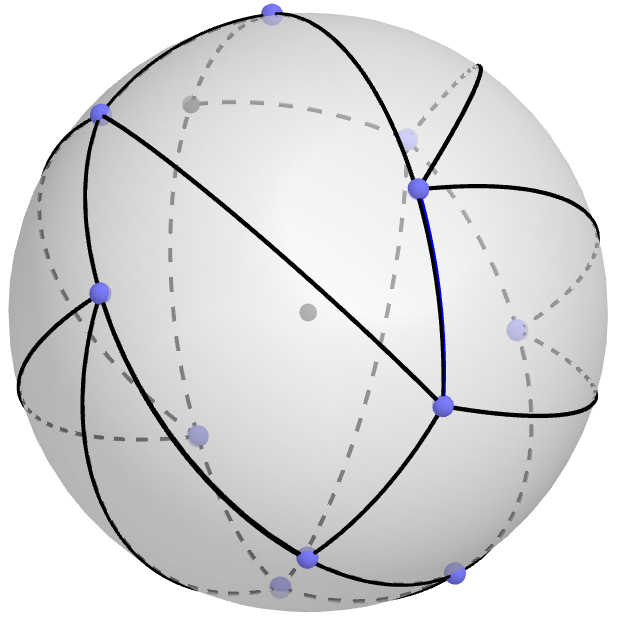}\hspace*{0.1cm}
	\includegraphics[width=0.133\textwidth]{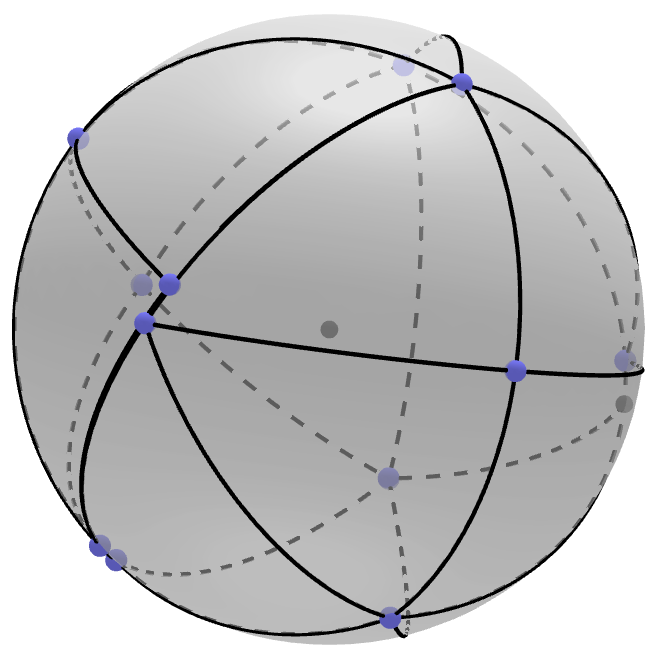}\hspace*{0.08cm}
	\includegraphics[width=0.133\textwidth]{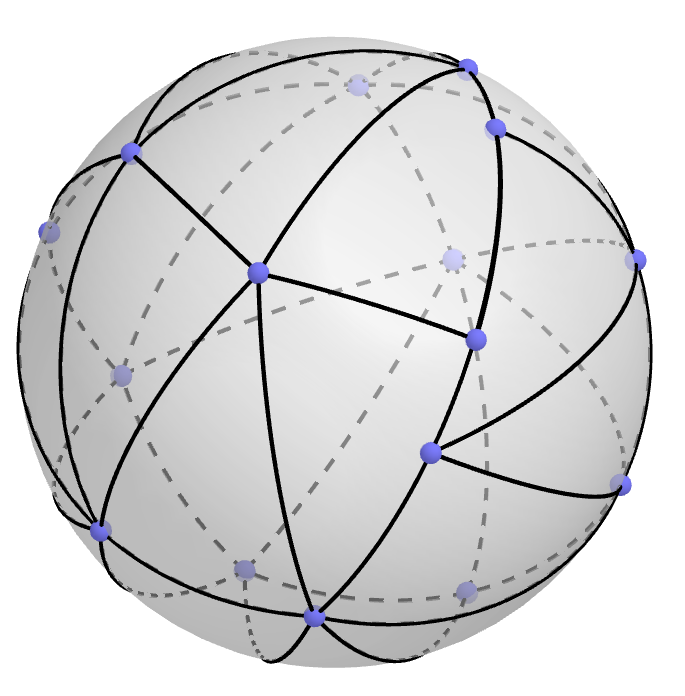}\hspace*{0.06cm}
	\includegraphics[width=0.133\textwidth]{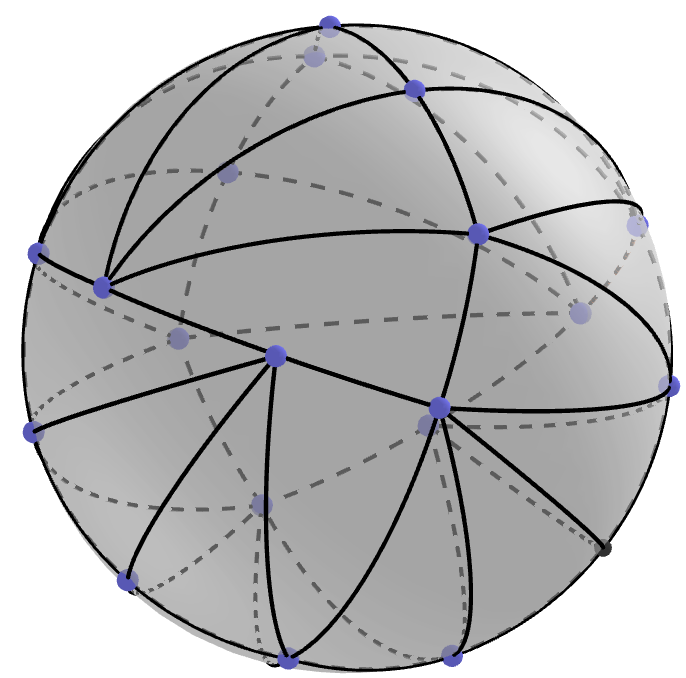}\hspace*{0.1cm}
	\includegraphics[width=0.133\textwidth]{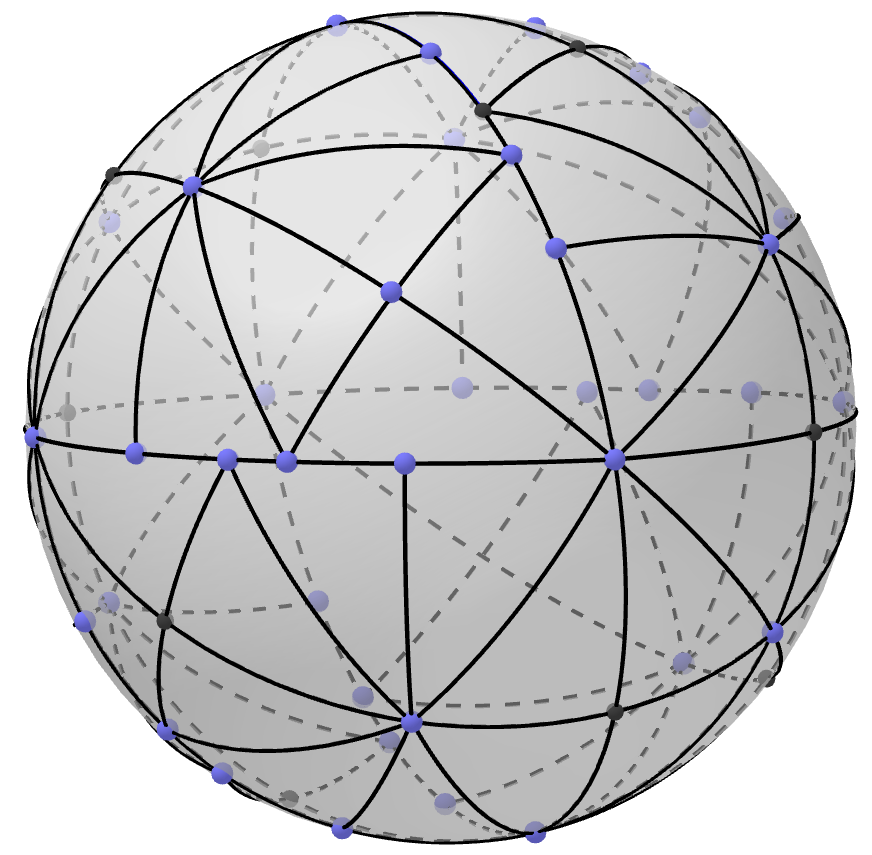}\hspace*{0.1cm}
	\includegraphics[width=0.133\textwidth]{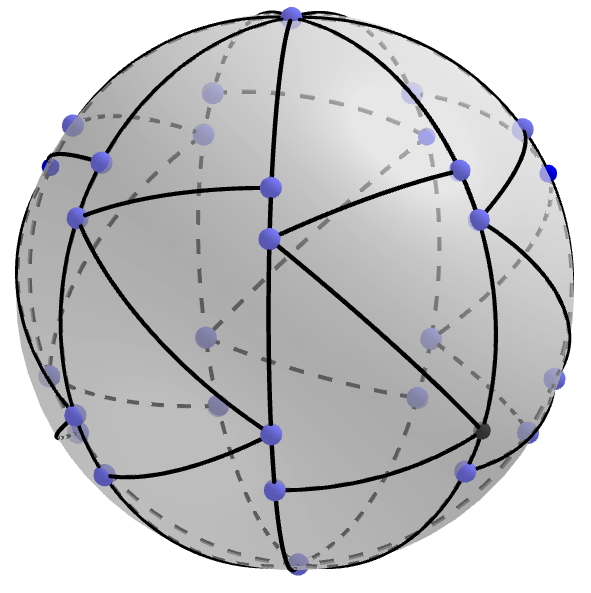}
	\caption{Non-side-to-side tilings by congruent right triangles.}
	\label{right figure}
\end{figure}

Seven more right triangle monotiles were found by Dawson and Doyle \cite{dawson2006,dawson2006-1,dawson2007}. The 3D pictures for some tilings are shown in Figure \ref{right figure}. The exact geometric data and the extended edges are given in Table \ref{right}.

\begin{table}[htp]
	\centering
	\caption{Seven more right triangle types and their tilings by Dawson and Doyle.}
	\label{right}
	\begin{tabular}{|c|c|c|}
		\hline
		$f$ & right triangle monotiles  & AVC and extended edges of some tiling\\
		\hline
		&$(\frac{3}{5}\pi,\frac{1}{2}\pi,\frac{3}{10}\pi)$&\\
		10&$a=\pi-\arccos(\frac{\frac{\sqrt{5}}{4}-\frac{1}{4}}{\frac{\sqrt{5}}{4}+\frac{1}{4}})\approx1.9627$&$\{4\aaa^2\bbb\ccc,2\aaa\bbb\ccc^3;\,2\bbb^2\}$\\
		&$b=\pi-\arccos(\frac{(\frac{\sqrt{5}}{4}-\frac{1}{4}) \sqrt{5-\sqrt{5}}}{\sqrt{5+\sqrt{5}}\, (\frac{\sqrt{5}}{4}+\frac{1}{4})})\approx 1.8091$&$b=2c$\\
		&$c=\arccos(\frac{\sqrt{5-\sqrt{5}}}{\sqrt{5+\sqrt{5}}})\approx0.9046$&\\
		
		\hline
		&$(\frac{7}{12}\pi,\frac{1}{2}\pi,\frac{1}{4}\pi)$&\\
		12&$a=\pi-\arccos(\frac{\sqrt{3}}{2}-\frac{1}{2})\approx1.9455$&$\{4\aaa^3\ccc;\,4\bbb^2,4\bbb\ccc^2\}$\\
		&$b=\pi-\arccos(\frac{\sqrt{3}-1}{1+\sqrt{3}})\approx1.8421$&$2a+c=c+2a$\\
		&$c=\arccos(\frac{2}{1+\sqrt{3}})\approx0.7495$&\\
		\hline
		
		&$(\frac{1}{2}\pi,\frac{5}{12}\pi,\frac{1}{3}\pi)$&\\
		16&$a=\arccos(\frac{(\sqrt{3}-1) \sqrt{3}}{3 (1+\sqrt{3})})\approx1.4155$&$\{2\aaa^4,\aaa^2\ccc^3,4\bbb^4\ccc;\,3\aaa^2,3\ccc^3\}$\\
		&$b=\arccos(\frac{\sqrt{2}\, (\sqrt{3}-1) \sqrt{3}}{6})\approx1.2673$&$2b=a+c$\\
		&$c=\arccos(\frac{\sqrt{2}}{1+\sqrt{3}})\approx 1.0267$&\\
		
		\hline
		&$(\frac{1}{2}\pi,\frac{5}{12}\pi,\frac{1}{4}\pi)$&\\
		24&$a=\arccos(\frac{\sqrt{3}-1}{1+\sqrt{3}})\approx1.2995$&$\{2\aaa^4,6\aaa\bbb^3\ccc,2\aaa^2\ccc^4,2\bbb^3\ccc^3;\,2\aaa\ccc^2,2\aaa^2\}$ \\
		&$b=\arccos(\frac{\sqrt{3}}{2}-\frac{1}{2})\approx1.1961$&$2a+c=c+2a$\\
		&$c=\arccos(\frac{2}{1+\sqrt{3}})\approx0.7495$&\\
		\hline
		
		&$(\frac{1}{2}\pi,\frac{7}{16}\pi,\frac{3}{16}\pi)$&\\
		32&$a=\arccos(\frac{\cos(\frac{7 \pi}{16}) \cos(\frac{3 \pi}{16})}{\sin(\frac{7 \pi}{16}) \sin(\frac{3 \pi}{16})})\approx1.2685$&$\{4\aaa^4,4\aaa^2\bbb\ccc^3,8\aaa\bbb^3\ccc;\,4\bbb\ccc^3\}$\\
		&$b=\arccos(\frac{\cos(\frac{7 \pi}{16})}{\sin(\frac{3 \pi}{16})})\approx1.2120
		$&$a+c=c+a$\\
		&$c=\arccos(\frac{\cos(\frac{3 \pi}{16})}{\sin(\frac{7 \pi}{16})})\approx0.5591$&\\
		\hline
		&$(\frac{1}{2}\pi,\frac{1}{3}\pi,\frac{2}{9}\pi)$&\\
		72&$a=\arccos(\frac{\cos(\frac{2 \pi}{9}) \sqrt{3}}{3 \sin(\frac{2 \pi}{9})})\approx0.8120$&$\{12\aaa^4,12\bbb^2\ccc^6,2\bbb^6;\,12\aaa^2,12\bbb^3\}$\\
		&$b=\arccos(\frac{1}{2 \sin(\frac{2 \pi}{9})})\approx0.6795$&$a+2c=2c+a$\\
		&$c=\arccos(\frac{2 \cos(\frac{2 \pi}{9}) \sqrt{3}}{3})\approx 0.4853$&\\
		\hline
		&$(\frac{1}{2}\pi,(\frac{1}{2}-\frac{1}{n})\pi,\frac{2}{n}\pi)^*$&\\
		4n&$a=\arccos\! \left(\frac{\cos\! \left(\frac{\left(f-8\right) \pi}{2 f}\right) \cos\! \left(\frac{8 \pi}{f}\right)}{\sin\! \left(\frac{\left(f-8\right) \pi}{2 f}\right) \sin\! \left(\frac{8 \pi}{f}\right)}\right)$&$\{2\ccc^n;\,2n\aaa^2,2n\bbb^2\ccc\}$\\
		&$b=\arccos\! \left(\frac{\cos\! \left(\frac{\left(f-8\right) \pi}{2 f}\right)}{\sin\! \left(\frac{8 \pi}{f}\right)}\right)$&$a+2b=2b+a$\\
		&$c=\arccos\! \left(\frac{\cos\! \left(\frac{8 \pi}{f}\right)}{\sin\! \left(\frac{\left(f-8\right) \pi}{2 f}\right)}\right)$&\\
		\hline		
	\end{tabular}
\end{table}

\subsubsection*{5. A sequence of triangle monotiles $(\frac{2k-1}{3k},\frac{1}{3},\frac{1}{k})$ $(k\geq 4)$ from quadrilateral tilings}

        In \cite{lw2}, the quadrilateral monotiles with corners ($\frac{2}{f'},\frac{4f'-4}{3f'},\frac{4}{f'},\frac{2f'-2}{3f'}$) ($f'$ represents the number of quadrilateral tiles) and three equal sides admit two-layer earth map tilings, and various modifications for $f' \equiv 4 \mod 6$. Such a quadrilateral can be subdivided into three congruent triangles as shown in Figure \ref{4bx}, with three angles   $\aaa=\frac{2f-6}{3f},\bbb=\frac{1}{3},\ccc=\frac{6}{f}$ ($f=3f'$ represents the number of triangle tiles). For example, the quadrilateral tilings in the first row of Figure \ref{4bxx}, after the subdivisions, become triangular tilings in the second row, with the AVC  $\{10\aaa^3\ccc,2\ccc^{10};\,10\bbb^3\}$,  $\{8\aaa^3\ccc,2\aaa^2\ccc^{4},2\aaa\ccc^7;\,10\bbb^3\}$,  $\{6\aaa^3\ccc,6\aaa^2\ccc^4;\,10\bbb^3\}$, $\{6\aaa^3\ccc,6\aaa^2\ccc^4,2\bbb^6;\,6\bbb^3\}$ respectively. The extended edges of these tilings are $a+c=c+a$.
	  \begin{figure}[H]
		\centering
		\begin{tikzpicture}[>=latex,scale=0.5]
			
	       \begin{scope}[xshift=9cm,yshift=-1cm]
	                  \draw (3.32, 1)--(8.04, 2.76)--(12.24, 7)--(11.66, 1);
	                  \draw[line width=1.5] (3.32, 1)--(11.66, 1);
	                  \draw[dashed] (8.04, 2.76)--(9.4,1)--(12.24, 7);
	                  \node at (4.92, 1.24){\footnotesize $\ccc$};\node at (11.2, 5.46){\footnotesize $\ccc$};\node at (11.76, 5.16){\footnotesize $\ccc$}; 
	                  \node at (11.28, 1.4){\footnotesize $\aaa$};\node at (8.54, 2.72){\footnotesize $\aaa$};\node at (7.92, 2.4){\footnotesize $\aaa$};
	                  \node at (7.8, 1.6){\footnotesize $\bbb=\frac13$};\node at(9.36, 1.68){\footnotesize $\frac13$}; \node at (10.26, 1.6){\footnotesize $\frac13$};
	                  \node at (5.44, 2.54){$b$};  \node at (9.5, 5.04){$b$}; \node at (12.58, 3.76){$b$}; \node at (8.04, 0.32){$a+c$};
	       \end{scope}
		\end{tikzpicture} 
\caption{The subdivision of a special quadrilateral monotile.}
\label{4bx}
\end{figure}
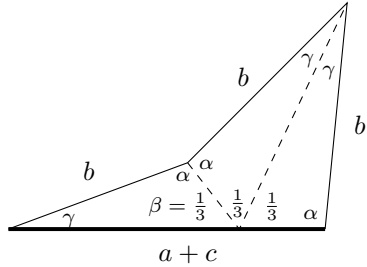
\begin{figure}[H]
	\centering
	\subfloat[10 quadrilaterals]{\includegraphics[width=0.21\textwidth]{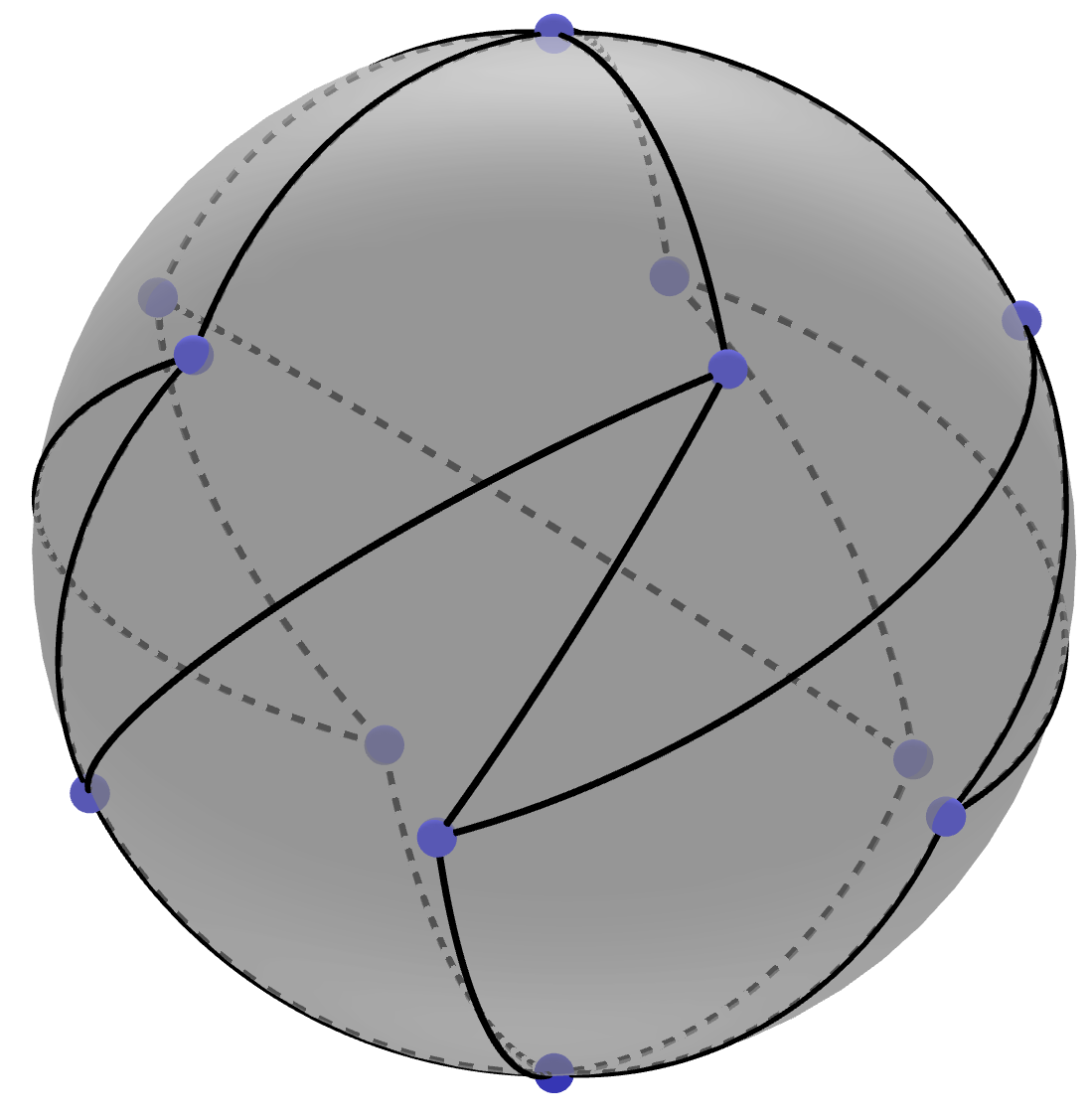}}\hspace*{0.3cm}
	\includegraphics[width=0.21\textwidth]{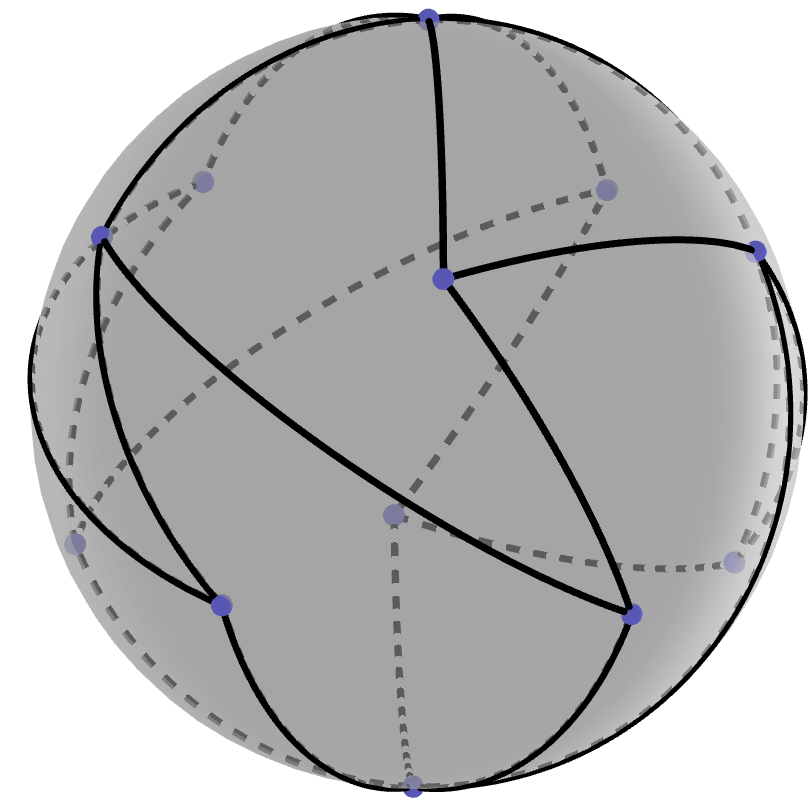}\hspace*{0.3cm}
	\includegraphics[width=0.21\textwidth]{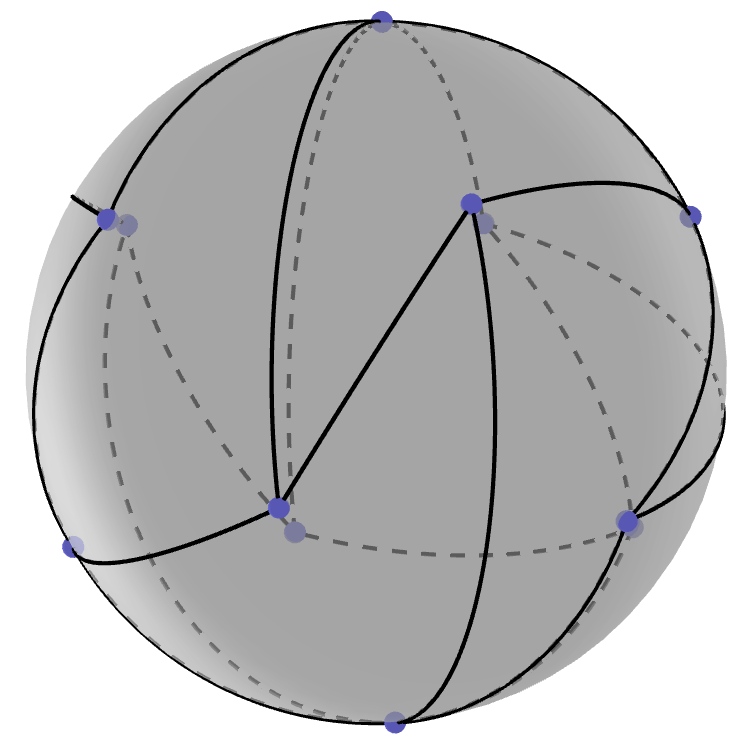}\hspace*{0.3cm}
	\includegraphics[width=0.21\textwidth]{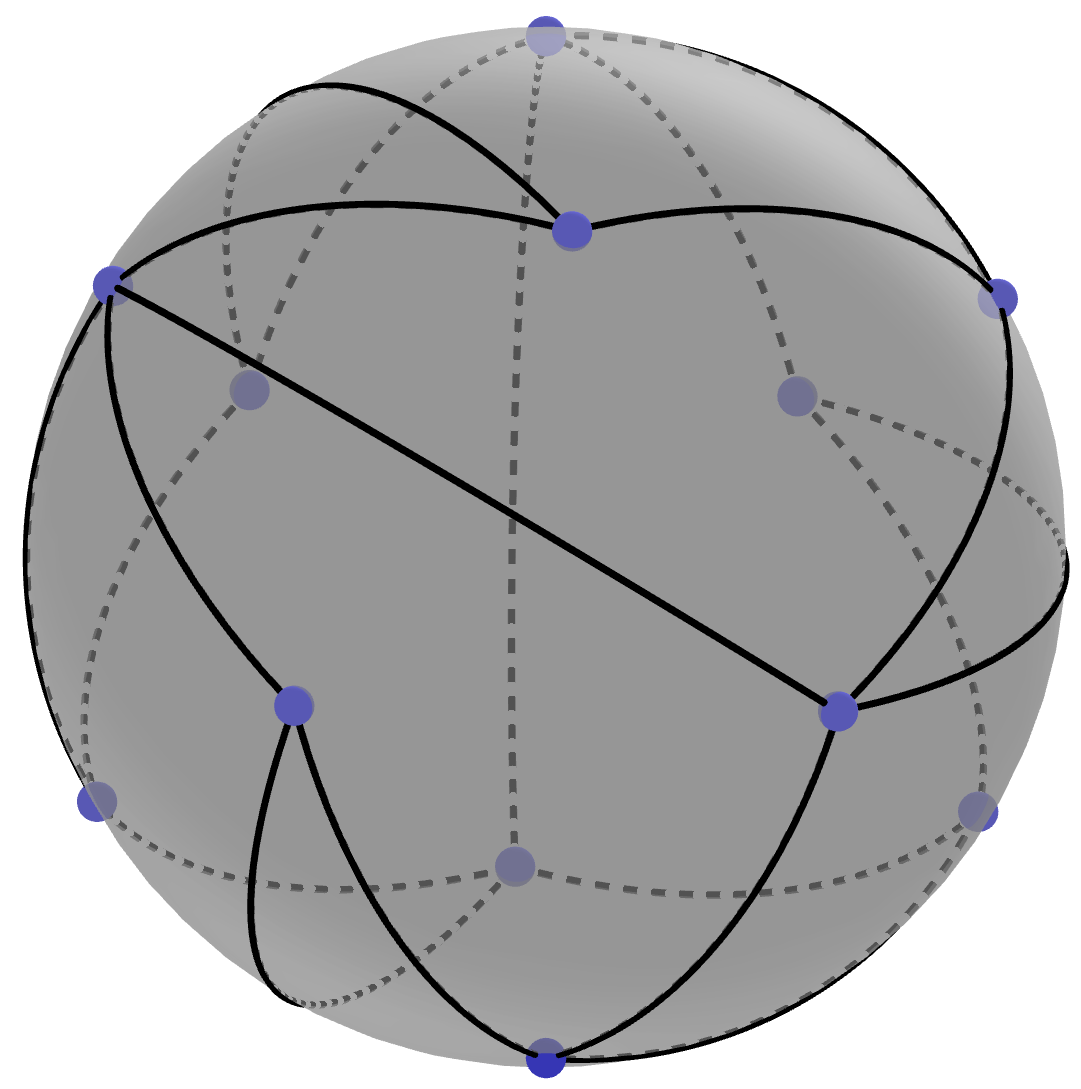}\\
	\subfloat[30 triangles]{\includegraphics[width=0.20\textwidth]{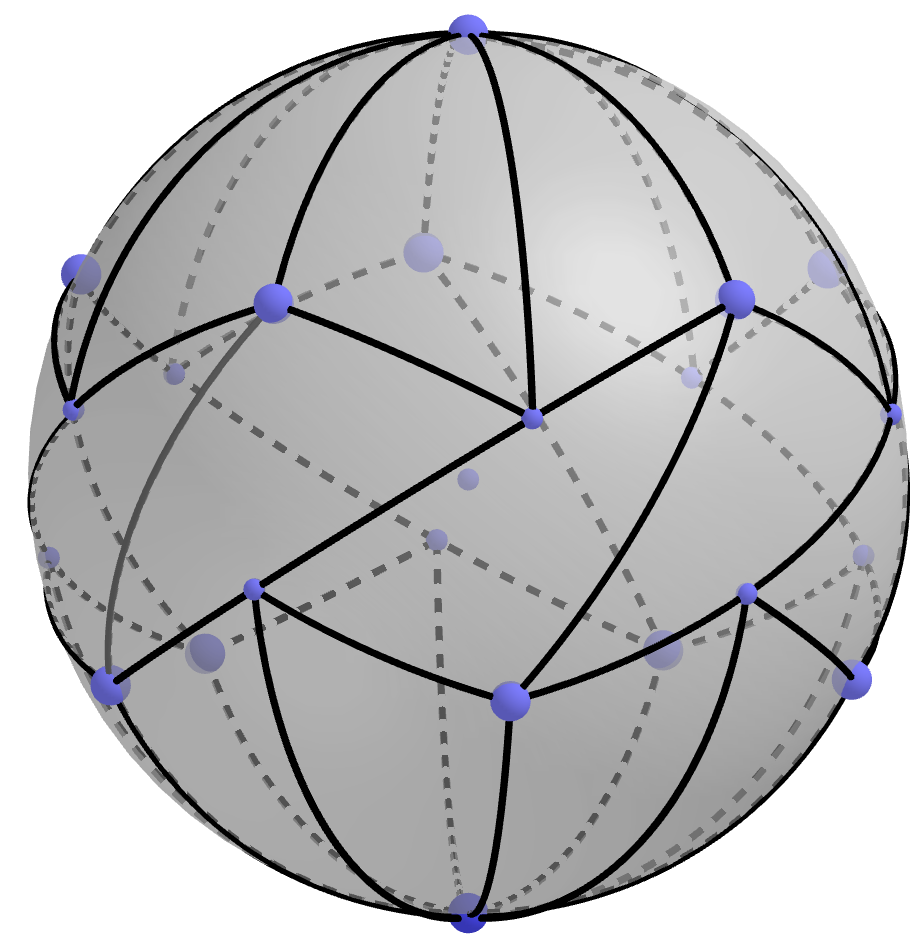}}\hspace*{0.3cm}
	\includegraphics[width=0.21\textwidth]{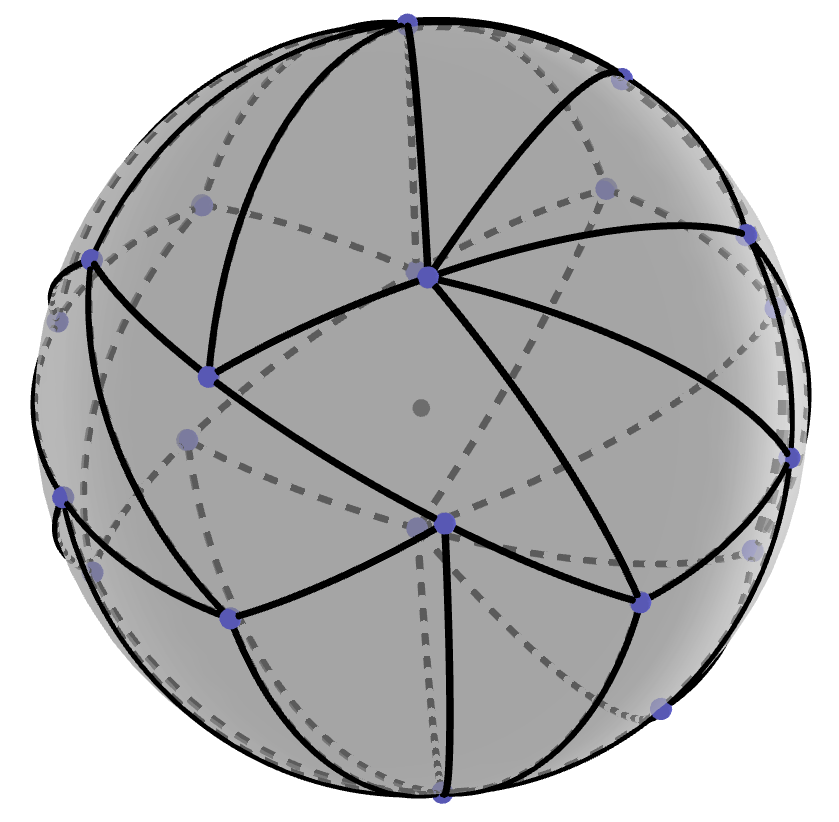}\hspace*{0.3cm}
	\includegraphics[width=0.22\textwidth]{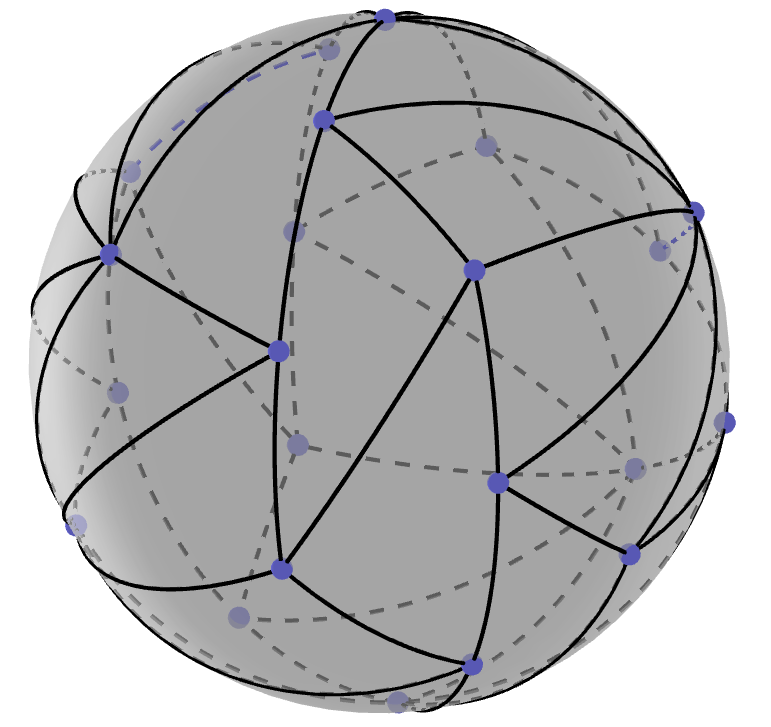}\hspace*{0.3cm}
	\includegraphics[width=0.20\textwidth]{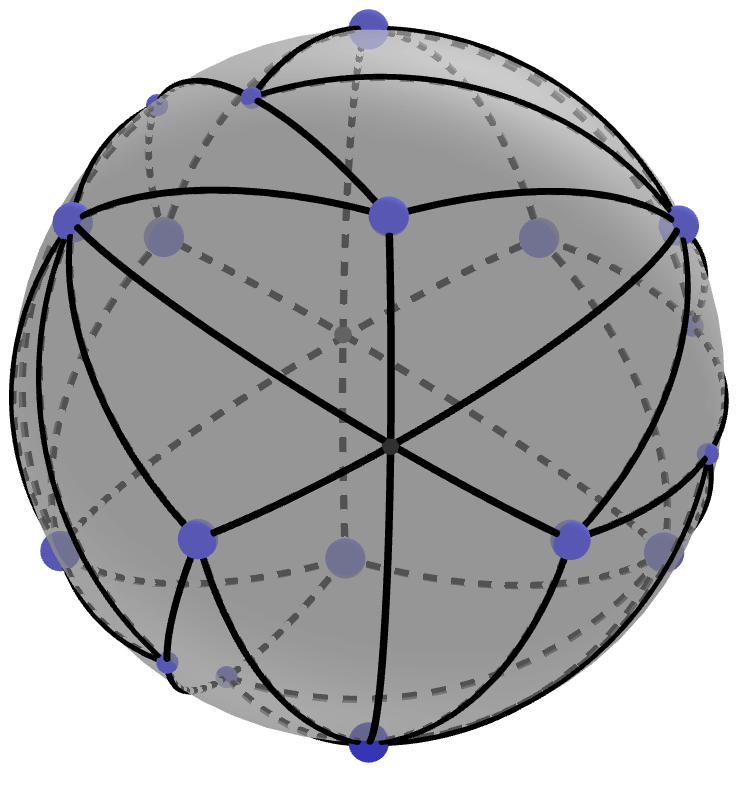}
	\caption{Non-side-to-side triangular tilings by subdividing  quadrilateral tilings.}
	\label{4bxx}

\end{figure}

In this paper, we extend Wang-Yan's systematic methods in  \cite{wy1,wy2}, such as statistics of vertex distribution and adjacent angle deduction, to handle non-side-to-side tilings, and get the following main result. 

\begin{theorem*}
All non-side-to-side tilings of the sphere by congruent triangles with any irrational angle are:

1. A sequence of one-parameter families of triangles admitting two-layer earth map tilings with $2k\ge 6$ triangles, and their rotation modifications for even $k$, as shown in Figure \ref{1.2}. 

2. A one-parameter family of triangles, each admitting a unique tiling with $8$ triangles, as shown in the first picture of Figure \ref{th-1}.

3. A unique triangle admitting a unique tiling with $16$ triangles, as  shown in the second picture of Figure \ref{th-1}.

All geometric data are given in Table \ref{data}.
	
\begin{figure}[H]
		\centering
		\includegraphics[width=0.2\textwidth]{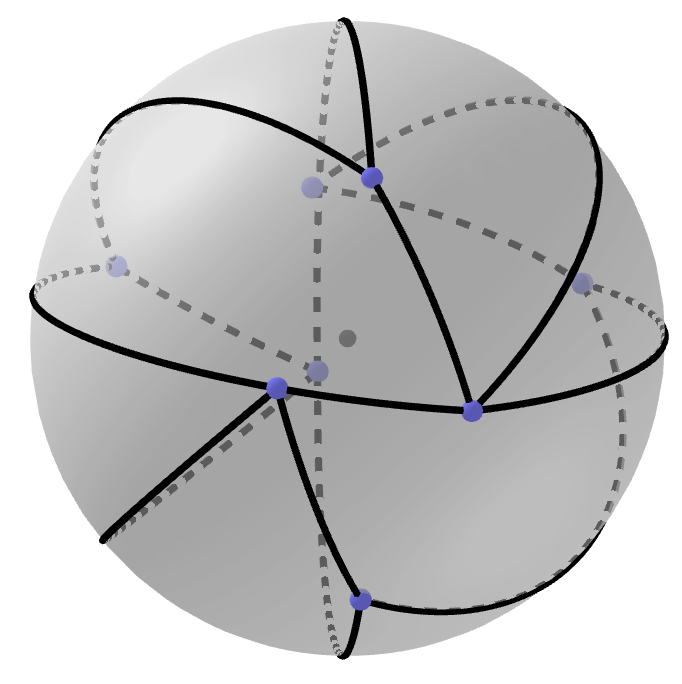}
		\hspace{1cm}  
		\includegraphics[width=0.203\textwidth]{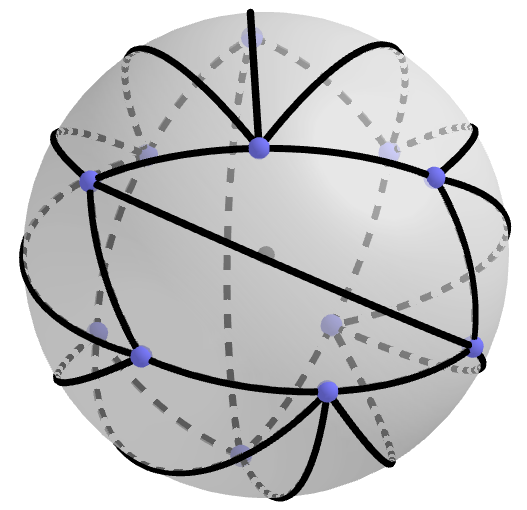}
		\caption{$T(4\aaa^2\ccc;4\bbb^2\ccc)$ and $T(8\aaa^2\bbb\ccc;4\bbb^2\ccc^2)$.}
		\label{th-1}
	\end{figure}
	\begin{table*}[htp]                    	
		\centering 
		\begin{tabular}{c|c}	 	
				angles $ (\aaa,\bbb,\ccc)$  & sides \\
				\hline 			 
				\multirow{3}{*}{$(\aaa,\pi-\aaa,\frac{4\pi}{f})$}&$a=\arccos\! \left (\frac{\cos\aaa(1-\cos\frac{4\pi}{f})}{\sin\frac{4\pi}{f}\sin\aaa}\right )$\\
				&$b=\pi-a$\\
				&$c=\arccos\left (\frac{\cos\frac{4\pi}{f}-\cos^2(\aaa)}{\sin^2(\aaa)}\right )$\\
				\hline
				\multirow{3}{*}{$(\pi-\frac{\ccc}{2},\frac{\pi}{2}-\frac{\ccc}{2},\ccc)$, $\ccc\in(\frac{\pi}{4},\frac{\pi}{2})$}
				&$a=\arccos\! \left(\frac{\tan\frac{\ccc}{2}\cos \ccc-1}{ \sin\ccc}\right)$\\
				\multirow{3}{*}{isosceles triangle when $\ccc=\frac\pi3$}&$b=\arccos\! \left(\frac{1-\cot\frac{\ccc}{2}\cos \ccc}{ \sin\ccc}\right)$\\
				&$c=\arccos(1-2\cot\ccc)$\\
				\hline
				\multirow{1}{*}{$(\frac{3}{4}\pi,\arctan\sqrt{2},\frac{\pi}{2}-\bbb)$}&$a=\frac{2 \pi}{3},b=\frac{\pi}{2},c=\frac{\pi}{4}$\\
				\hline
				\end{tabular}
				
				\caption{The exact geometric data for the triangle monotiles in the theorem.}
				\label{data}
			\end{table*}    

\label{th}
\end{theorem*}


{\bf Acknowledgement}\ \  The authors would like to thank Professor Min Yan at HKUST and a graduate student Yixi Liao at ZJNU for many helpful discussions. The results of this paper were reported on the 11TH Shanghai Conference on Combinatorics (11SHCC) in May 24 - 28, 2024. The authors would like to thank Shanghai Jiao Tong University and the organizers for great hospitality.

\section{Basic Facts}

\subsection*{Corner, side of a polygon and vertex, edge of a tiling}

\begin{figure}[H]
	\centering
	\begin{tikzpicture}[>=latex,scale=0.5]
		\draw (6, 5)--(9, 5)--(10, 5)--(6, 3)--(6, 5)--(3, 3)--(6, 3)
		(3.59,3.39)--(3.86, 6.42)--(6, 5)--(7, 7)--(7.46, 7.94)--(3.86, 6.42)
		(9, 5)--(7, 7)
		(7, 7)--(12, 5)--(10, 5);
		\fill 
		(6, 5) circle (0.1)
		(9, 5) circle (0.1);
		
		\node at (6.38, 4.66) {\small A};
		\node at (9.26, 5.36) {\small B};

	\end{tikzpicture} 
	\caption{Full vertex $A$ and half vertex $B$.}
	\label{vh}
\end{figure}
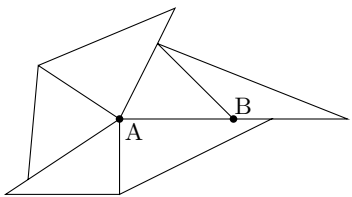

For the basics of tilings, we adopt the definitions and conventions in the books \cite{gs} and \cite{ac1}. We emphasize that corners and sides refer to polygons or tiles, and vertices and edges refer to the tiling. If the vertices and edges coincide with the corners and sides, then we say the tiling is \emph{side-to-side}\footnote{This is usually called {\em edge-to-edge} in the literature.}. If a tiling is not side-to-side, then a vertex may lie in the interior of a side, such as the vertex $B$ in Figure \ref{vh}. Then we call the vertex a {\em half vertex}. Otherwise the vertex is a {\em full vertex}, such as the vertex $A$ in Figure \ref{vh}. A vertex is a full vertex if and only if it is the corner of all the tiles at the vertex.

The {\em degree} of a vertex is the number of corners at the vertex. For example, the full vertex $A$ in Figure \ref{vh} has degree $5$, and the half vertex $B$ has degree $2$. In this paper, we only consider tilings in which all sides are straight.

\subsection*{Stop Vertex and Extended Edges}

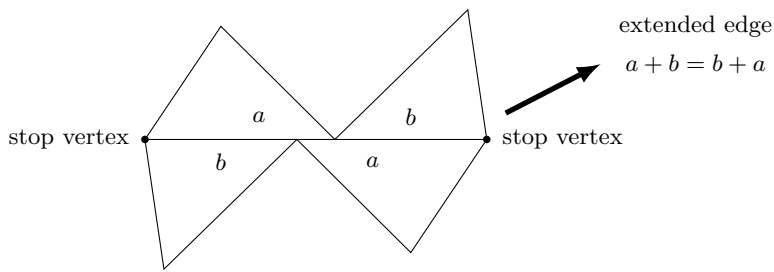
\begin{figure}[H]
 \centering
 \begin{tikzpicture}[>=latex,scale=0.5]
  \draw (0,0)--(2,3)--(5,0)--(8.5,3.44)--(9,0)--(7,-3)--(4,0)--(0.5,-3.44)--(0,0)--(9,0);
  \fill 
  (0, 0) circle (0.1)
  (9, 0) circle (0.1);
  
   \node at (3, 0.6) {\small $a$};
   \node at (2, -0.6) {\small $b$};
   \node at (7, 0.6) {\small $b$};
   \node at (6, -0.6) {\small $a$};
  \node at (-2, 0) {\small stop vertex};
  \node at (11, 0) {\small stop vertex};
  \node at (14.5, 3) {\small extended edge};

  \draw[line width=2pt, ->](9.5,0.7)--(12,2);
  \node at (14.5, 2) {\small $a+b=b+a$};
 \end{tikzpicture} 
 \caption{Mismatched tiles, an extended edge, and two stop vertices.}
 \label{extended}
\end{figure}

The appearance of a half vertex means two adjacent tiles are not in the side-to-side position. We describe this situation as the tiles being $unmatched$ or $mismatched$. In this case, the side can be extended to a straight line consisting of many sides in the tiling, which we call $extended$ $edges$. In each of the two directions of the extended edges, the edges either extend forever or contain a vertex that is the corner of both adjacent tiles on the two sides of the edges. This vertex is called a $stop$ $vertex$ (along the extended edges).   In an extended edge, if there is only one stop vertex or no stop vertex, then it must be a great circle of length $2\pi$, like the seventh picture in Figure \ref{subdivision}. An extended edge between two stop vertices is denoted by an equality of ordered sides, as shown in Figure \ref{extended}. 

\subsection*{Adjacent Angle Deduction}
The very useful tool \textit{adjacent angle deduction} (abbreviated as AAD) has been introduced in \cite[Section 2.5]{wy1}. The following is \cite[Lemma 10]{wy1}. We give a quick review here using Figure \ref{AAD}. Let “$\thin$" denote an side, $\thin^\bbb\aaa^\ccc\thin$ denote a triangle. Then we indicate the arrangement of angles and sides by denoting the vertices as $\thin^\ccc\bbb^\aaa\thin^\ccc\aaa^\bbb\thin^\bbb\aaa^\ccc\thin^\aaa\ccc^\bbb\thin$. The notation can be reversed, such as $\thin^\ccc\bbb^\aaa\thin^\ccc\aaa^\bbb\thin^\bbb\aaa^\ccc\thin^\aaa\ccc^\bbb\thin$=$\thin^\bbb\ccc^\aaa\thin^\ccc\aaa^\bbb\thin^\bbb\aaa^\ccc\thin^\aaa\bbb^\ccc\thin$; and it can be rotated, such as $\thin^\ccc\bbb^\aaa\thin^\ccc\aaa^\bbb\thin^\bbb\aaa^\ccc\thin^\aaa\ccc^\bbb\thin$= $\thin^\ccc\aaa^\bbb\thin^\bbb\aaa^\ccc\thin^\aaa\ccc^\bbb\thin^\ccc\bbb^\aaa\thin$.  If we assume $a>b>c$, then $\aaa^\ccc\thin^\aaa\bbb$ implies a half vertex at $\aaa$.

\begin{figure}[H]
	\centering
	\begin{tikzpicture}[>=latex,scale=0.4]
		\draw (3, 6)--(7,6)--(6.71,2.3)--(4.9, 4.14)--(7,6)
		(3, 6)--(3, 2)--(7, 2)--(6.71,2.3)--(3, 6)
		(3, 2)--(4.9, 4.14);
		
		\node at (6.48, 3) {\small $\aaa$};
		\node at (5.44, 4.12) {\small $\bbb$};
		\node at (6.68, 5.28) {\small $\ccc$};
		
		\node at (4.88, 3.66) {\small $\aaa$};
		\node at (3.88, 2.42) {\small $\bbb$};
		\node at (6.1, 2.36) {\small $\ccc$};
		
		\node at (4.44, 4.14){\small $\aaa$};
		\node at (3.4, 3) {\small $\bbb$};
		\node at (3.4, 5.12) {\small $\ccc$};
		
		\node at (3.94, 5.68){\small $\aaa$};
		\node at (6.12, 5.68) {\small $\bbb$};
		\node at (4.98, 4.64) {\small $\ccc$};
		
	\end{tikzpicture} 
	\caption{ Adjacent Angle Deduction of $\thin\bbb\thin\aaa\thin\aaa\thin\ccc\thin$. }
	\label{AAD}
\end{figure}

\begin{lemma}\label{angle sum}
	If all tiles in a tiling of the sphere by $f$ triangles have the same three corners $\aaa$, $\bbb$, $\ccc$, then
	\begin{center}
		$\aaa+\bbb+\ccc =(1+\frac{4}{f})$.
	\end{center}
	
\end{lemma}

\begin{lemma}[Balance Lemma]\label{balance lemma}
	If all tiles in a tiling of the sphere by $f$ triangles have the same three corners $\aaa$, $\bbb$, $\ccc$, then a vertex (either full or half) having more $\aaa$ than $\bbb$ (or $\ccc$) implies the existence of some other vertex having more $\bbb$ (or $\ccc$) than $\aaa$. 
\end{lemma}

\begin{lemma}\label{degen}
	In a non-side-to-side tiling of the sphere by congruent triangles, if any angle is $\pi$, then the triangle monotile  degenerates to a 2-gon or 1-gon (a hemisphere), and the tiling is an earth map tiling or its rotation modification. 
\end{lemma}

\begin{proof}
		Assume $\aaa=\pi$, then we have $b \cup c$ and $a$ are two different geodesics (great arcs) between $\bbb$ and $\ccc$. Either three edges form a great circle: $a+b+c=2\pi$; or $\bbb$ is the antipode of $\ccc$: $a=b+c=\pi$.
		
		In the first case, we have $\bbb=\ccc=\pi$ and  $f=2$. Note that $b$, $c$ are not equal in general. The tiles are two hemispheres which may rotate along the equator to any position, as shown in the first picture of Figure \ref{aaa=1}. Each vertex is either a degree 2 full vertex or a degree 1 half vertex. 
		
		In the second case, the spherical triangle degenerates into a 2-gon with  $\bbb=\ccc=\frac{2\pi}{f}$. Note that $b$, $c$ are not equal in general. When $f$ is odd, each $\bbb$ must belong to a full vertex  $\bbb^k\ccc^{f-k}$, which induces the antipodal full vertex $\bbb^{f-k}\ccc^k$. The tiling is an earth map tiling with each tile being a standard timezone, and each $\aaa$-vertex is either a degree 2 full vertex or a degree 1 half vertex, as shown in the second picture of Figure \ref{aaa=1}. When $f$ is even, $\frac{f}{2}$ continuous timezones form a hemisphere, which may rotate along the equator to any position, as shown in the third picture (with a full vertex $\aaa\bbb^k\ccc^{\frac{f}{2}-k}$) and the fourth picture (with four half vertices $\bbb^k\ccc^{\frac{f}{2}-k}$, $\bbb^{\frac{f}{2}-k}\ccc^k$,$\bbb^l\ccc^{\frac{f}{2}-l}$, $\bbb^{\frac{f}{2}-l}\ccc^l$) of Figure \ref{aaa=1} respectively.	
		\begin{figure}[H]
			\centering
			\includegraphics[width=0.165\textwidth]{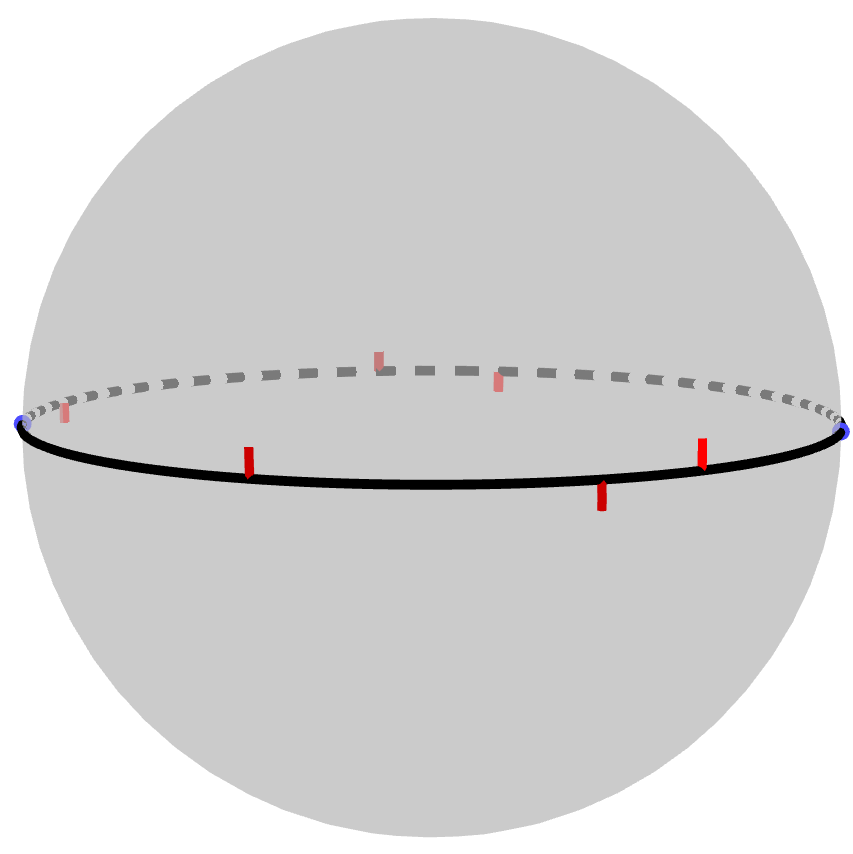}\hspace*{0.3cm}
			\includegraphics[width=0.154\textwidth]{xxaa6-1}\hspace*{0.3cm}
			\includegraphics[width=0.161\textwidth]{xxaa6-2}\hspace*{0.3cm}
			\includegraphics[width=0.156\textwidth]{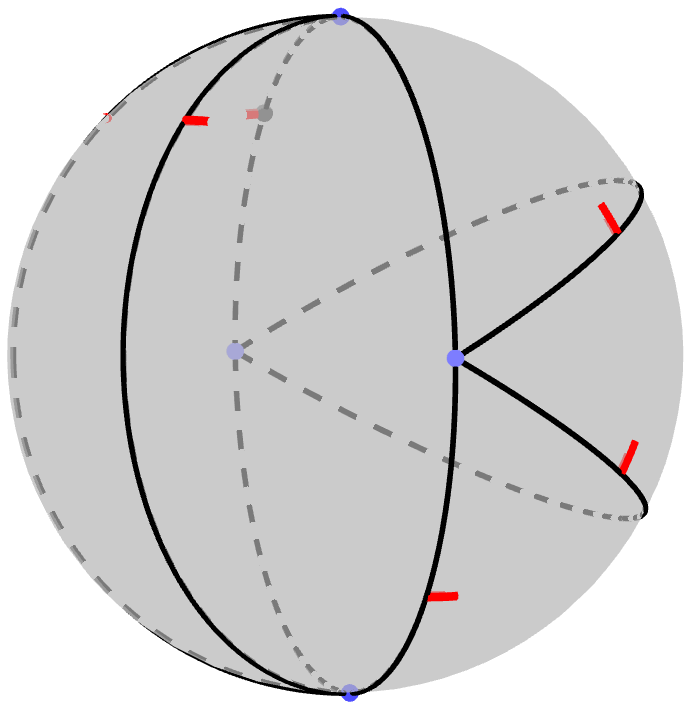}
			\caption{All possible non-side-to-side tilings by congruent degenerate triangles.}
			\label{aaa=1}
		\end{figure}

\end{proof}

In a non-side-to-side tiling of the sphere by congruent triangles, let $v,h,e,f$ be the numbers of full vertices, half vertices, edges, and tiles. Then we have
\begin{align}
	(v+h)-e+f &=2, \label{eq1} \\
	3f+h &=2e. \label{eq2}
\end{align}

The first equality is the Euler number of the sphere. For the second equality, we note that the number of edges in a tile is $3$ plus the number of half vertices on the sides of the tile. Since each half vertex has only one side, corresponding to only one tile, the sum $\sum$ of all these numbers is $3f+h$. On the other hand, since each edge is shared by exactly two tiles, the sum $\sum$ is $2e$. 

Canceling $e$ from the equations \eqref{eq1} and \eqref{eq2}, we get $f=2v+h-4$.

 Let $v_k$ be the number of full vertices of degree $k$ ($k\geq2$). Let $h_l$ be the number of half vertices of degree $l$ ($l\geq1$). Then we have
\begin{align*}
	v &=\sum_{k\ge 2}v_k=v_2+v_3+v_4+v_5+\cdots, \\
	h &=\sum_{l\ge 1}h_l=h_1+h_2+h_3+h_4+\cdots, \\
	3f &=\sum_{k\ge 2}kv_k+\sum_{l\ge 1}lh_l.
\end{align*}
The first two equalities follow directly from the definition. The last equality is due to two ways of counting the total number of corners in a tiling. Substituting the three equalities into $f=2v+h-4$, we get 
\begin{align}
	4v_2&+3v_3+2v_4+v_5+2h_1+h_2=12+(v_7+2v_8+\cdots)+(h_4+2h_5+\cdots), \label{eq3} \\
	3f&+2v_2+h_1=12+(2v_4+4v_5+6v_6+\cdots)+(h_2+3h_3+5h_4+\cdots). \label{eq4}
\end{align}

These formulas do not change if we replace any full vertices in $v_{2l}$ with two half vertices in $h_l$.

\begin{lemma}\label{convex}
	In a non-side-to-side tiling of the sphere by congruent triangles, the triangle is not concave.
\end{lemma}

\begin{proof}
	 Otherwise assume $\aaa>\pi$, which implies $\bbb,\ccc\neq\pi$ and $h_1=0$ by Lemma \ref{degen}. Then there is no $\aaa$ in any half vertex, at most one $\aaa$ in any full vertex. So each of the $f$  $\aaa$'s must appear in a different full vertex, which forces $v\ge f$. By $0<h=f-2v+4 \le 4-f$ and  \eqref{eq4}, we have $f=2$ or $3$. 
	  
	  If $f=2$, then $2v_2\ge 6$ by \eqref{eq4}, which implies that all $6$ angles of the two tiles are counted in $v_2$. There is no half vertex, a contradiction. 
	 
	 If $f=3$, then we have $h=1$, $v=3$,  $\aaa+\bbb+\ccc=\frac{7}{3}\pi$. Then \eqref{eq4} implies either $v_2=2,v_3=1,h_2=1$ or $v_2=3, h_3=1$. The AVC must be $\{2\aaa\bbb, \aaa\ccc^2; \bbb\ccc\}$ or $\{3\aaa\bbb; \ccc^3\}$. In the first case, we get $\aaa=\frac{4}{3}$, $\bbb=\frac{2}{3}$, $\ccc=\frac{1}{3}$, and this triangle does not exist. For the second case, we have AAD $\ccc^3=\thin^{\aaa}\ccc^{\bbb}\thin^{\aaa}\ccc^{\bbb}\thin^{\aaa}\ccc^{\bbb}\thin$.  Then we get $a=b$, and this triangle does not exist. 

	 	 \end{proof}

Since both degenerate and concave cases are clarified in Lemma \ref{degen} and Lemma \ref{convex}, it is enough to consider convex triangles henceforth. 

\begin{lemma}[Ref. \cite{ac1}]\label{Dawosn}
The corner values and side lengths of a convex spherical triangle satisfies  $0<\aaa,\bbb,\ccc<1$, $0<a,b,c<1$ and
	\[
	\aaa+\bbb<1+\ccc,\quad \bbb+\ccc<1+\aaa,\quad \aaa+\ccc<1+\bbb;
    \]
	\[
	a+b>c,\quad a+c>b,\quad b+c>a.
	\]

\end{lemma}

Dawson has completely classified the tilings of the sphere by congruent equilateral or isosceles triangles in \cite{dawson2001,dawson2003}. Therefore, we only need to discuss convex and scalene triangle, and always assume $\aaa>\bbb>\ccc$, $a>b>c$. Throughout this paper, an $abc$-tiling is always a  tiling of the sphere by congruent convex and scalene triangles in Figure 1. Then all full vertices have degree $\geq$ 3 and all half vertices have degree $\geq$ 2.

\begin{lemma}\label{345}
	In a non-side-to-side $abc$-tiling, $f>4$ and a full vertex of degree 3,4,5 or a half vertex of degree 2 in Table \ref{vertex} must appear.
\end{lemma}
\begin{proof}
	By \eqref{eq3}, we have $3v_3+2v_4+v_5+h_2\geq12$. This implies that there must exist a full vertex of degree 3,4,5 or a half vertex of degree 2. In a non-side-to-side tiling, we have $h>0$. Therefore, we get $f>4$ by \eqref{eq4}.

	By $f>4$, we get $\aaa+\bbb+\ccc=1+\frac{4}{f}<2$. This implies $\aaa\bbb\ccc$, $\aaa\ccc^2$, $\bbb^2\ccc$, $\bbb\ccc^2$, $\ccc^3$ are not vertices. For all other degree 3 vertices, we have $\aaa^3,\bbb^3,\aaa^2\bbb,\aaa\bbb^2,\aaa^2\ccc$.
	Together with all possibilities of degree 4,5 full vertices and degree 2 half vertices, we get Table \ref{vertex}
\end{proof}
\begin{table}[H]
	\centering
	\caption{The classification of degree 3,4,5 full vertices and degree 2 half vertices.}
	\renewcommand{\arraystretch}{1.5} 
	
	\begin{tabular}{c|c}
		\hline
		$v_3$ & $\aaa^3,\bbb^3,\aaa^2\bbb,\aaa\bbb^2,\aaa^2\ccc$ \\
		
		$v_4$ & $\aaa^3\bbb, \aaa^3\ccc, \aaa\bbb^3, \aaa\ccc^3, \bbb^3\ccc, \bbb\ccc^3, \aaa^2\bbb\ccc, \aaa\bbb^2\ccc, \aaa\bbb\ccc^2,\aaa^4,\bbb^4,\ccc^4,\aaa^2\bbb^2,\aaa^2\ccc^2,\bbb^2\ccc^2$ \\
		
		$v_5$ & $\aaa ^5, \bbb^5, \ccc^5, \aaa^4\bbb, \aaa^3\bbb^2, \aaa^2\bbb^3, \aaa\bbb^4, \aaa^4\ccc, \aaa^3\ccc^2, \aaa^2\ccc^3, \aaa\ccc^4, \bbb^4\ccc, \bbb^3\ccc^2, \bbb^2\ccc^3, \bbb\ccc^4,$ \\
		& $\aaa^3\bbb\ccc, \aaa^2\bbb^2\ccc, \aaa^2\bbb\ccc^2, \aaa\bbb^3\ccc, \aaa\bbb\ccc^3, \aaa\bbb^2\ccc^2$ \\
		
		$h_2$ & $\aaa^2,\bbb^2,\ccc^2,\aaa\bbb,\aaa\ccc,\bbb\ccc$ \\
		\hline
	\end{tabular}
\label{vertex}
\end{table}

 A vertex $\aaa^{n_1}\bbb^{n_2}\ccc^{n_3}$ can also be represented by its vector type $\boldsymbol{n}=(n_1,n_2,n_3)\in \mathbb{N}^3$. 

\begin{lemma}[Irrational Angle Lemma]\label{irrational}
	
	 In an $abc$-tiling with any irrational angle and $f$ tiles, if there is a vertex $\boldsymbol{m}\in \mathbb{N}^3$ and $\boldsymbol{m}$, $\boldsymbol{u}:=(1, 1 , 1 )$ are linearly independent, then all vertices $\bn$ belong to 
	\[
	 \mathbb{N}^3\cap\mathrm{Line}\left\{\boldsymbol{m},\frac{2f}{f+4}\boldsymbol{u}\right\}. 
	\] 
	In particular, the square matrix of $\boldsymbol{u},\boldsymbol{m},\boldsymbol{n}$ has determinant 0.
\end{lemma}
\begin{proof} 
	By Lemma 1 and the vertices $\boldsymbol{m}, \boldsymbol{n}$, the angles satisfy a linear system of equations with the augmented matrix 
	\[C=\left(\begin{matrix}
		1 & 1 & 1 &  1+\tfrac4f \\
		m_1 & m_2 & m_3 &  2\\
		n_1 & n_2 & n_3 &  2
	\end{matrix}\right).
	\]
	If $\bn$ is not a linear combination of $\boldsymbol{m}, \boldsymbol{u}$, then the system has a unique rational solution, a contradiction. So the last row of $C$ must be some linear combination of the other $2$ rows, which implies that $\bn$ lies in the line containing the two points $\boldsymbol{m},\frac{2f}{f+4}\boldsymbol{u}$. 
\end{proof}

\begin{lemma}[Matching Lemma] \label{mis_M}
In a non-side-to-side $abc$-tiling, if $\aaa$ never appears at any half vertex, then all possible extended edges are given in Table 6.

\begin{table}[H]
	\centering
	\caption{All possible extended edges if $\aaa$ is not in any half vertex.}
	\label{cuo}
\begin{tabular}{c|c||c|c}
	\hline
	& Extended Edge&  & Extended Edge \\
	\hline
	1& $b+ka=ka+b$ &  5&  $(k+1)a=c+ka+c$ \\
	\hline
	2& $c+ka=ka+c$ &  6& $c+ka+c=ka+b$ \\
	\hline
	3& $c+ka+b=b+ka+c$ &  7& $c+(k+1)a+c=b+ka+b$ \\
	\hline
	4& $b+ka+b=(k+1)a+c$ &  8& great circle \\
	\hline
\end{tabular}
\end{table}

\begin{figure}[htp]
	\centering
	\begin{tikzpicture}[>=latex,scale=1]
		
		
		\draw
		(0,0) -- (2,0)
		(0,-0.5) -- (0,0.5)
		(1,0) -- (1,0.5)
		(2,-0.5) -- (2,0.5);
		
		\fill 
		(0,0) circle (0.05)
		(2,0) circle (0.05);
		
		\node at (0.15,0.2) {\small $\alpha$};
		\node at (0.85,0.2) {\small $\beta$};
		\node at (0.5,0.2) {\small $c$};
		
		\node at (1.85,0.2) {\small $\alpha$};
		\node at (1.15,0.2) {\small $\beta$};
		\node at (1.5,0.2) {\small $c$};
		
		\node at (1,-0.2) {\small $a$};
		
		
		\begin{scope}[xshift=2.5cm]
			
			\draw
			(0,0) -- (2.5,0)
			(0,-0.5) -- (0,0.5)
			(1,0) -- ++(0,0.5)
			(2.5,0) -- ++(0,0.5)
			(2,0) -- ++(0,-0.5);
			
			\fill 
			(0,0) circle (0.05);
			
			\node at (0.15,0.2) {\small $\alpha$};
			\node at (0.85,0.2) {\small $\beta$};
			\node at (0.5,0.2) {\small $c$};
			\node at (2.35,0.2) {\small $\alpha$};
			\node at (1,-0.2) {\small $a$};
			
			\node at (1,-0.6) {\small no};
			
		\end{scope}

		
		\begin{scope}[xshift=5.5cm]
			
			\draw
			(0,0) -- (3,0)
			(0,-0.5) -- (0,0.5)
			(1,0) -- ++(0,0.5)
			(3,-0.5) -- (3,0.5)
			(2,0) -- ++(0,-0.5);
			
			\fill 
			(0,0) circle (0.05)
			(3,0) circle (0.05);
			
			\node at (0.15,0.2) {\small $\alpha$};
			\node at (0.85,0.2) {\small $\beta$};
			\node at (0.5,0.2) {\small $c$};
			\node at (2,0.2) {\small $a$};
			\node at (1,-0.2) {\small $a$};
			\node at (2.5,-0.2) {\small $c$};
			\node at (2.85,-0.2) {\small $\alpha$};
			\node at (2.15,-0.2) {\small $\beta$};
			
		\end{scope}
		
		
		\begin{scope}[xshift=9cm]
			
			\draw
			(0,0) -- (4,0)
			(0,-0.5) -- (0,0.5)
			(1,0) -- ++(0,0.5)
			(3,0) -- ++(0,0.5)
			(2,0) -- ++(0,-0.5)
			(4,0.5) -- (4,-0.5);
			
			\fill 
			(0,0) circle (0.05)
			(4,0) circle (0.05);
			
			\node at (0.15,0.2) {\small $\alpha$};
			\node at (0.85,0.2) {\small $\beta$};
			\node at (0.5,0.2) {\small $c$};
			\node at (2,0.2) {\small $a_1$};
			\node at (1,-0.2) {\small $a$};
			\node at (3,-0.2) {\small $b$};
			\node at (3.85,-0.2) {\small $\alpha$};
			\node at (2.15,-0.2) {\small $\gamma$};
			\node at (3.85,0.2) {\small $\alpha$};
			\node at (3.15,0.2) {\small $\beta$};
			\node at (3.5,0.2) {\small $c$};

		\end{scope}
		
	\end{tikzpicture}
	\caption{Matching of $a+\cdots$ and $c+\cdots$ along an extended line}
	\label{ac_matching}
\end{figure}

\begin{figure}[htp]
	\centering
	\begin{tikzpicture}[>=latex,scale=1]
		
		\draw
		(0,0) -- (3,0)
		(0,-0.5) -- (0,0.5)
		(2,0) -- (2,0.5)
		(2.5,0) -- (2.5,-0.5)
		(3,0) -- (3,0.5);
		
		\fill 
		(0,0) circle (0.05)
		(3,0) circle (0.05);
		
		\node at (0.15,0.2) {\small $\alpha$};
		\node at (1.85,0.2) {\small $\ccc$};
		\node at (1,0.2) {\small $b$};
		
		\node at (2.15,0.2) {\small $\bbb$};
		\node at (2.85,0.2) {\small $\aaa$};
		\node at (2.5,0.2) {\small $c$};
		
		\node at (1,-0.2) {\small $a$};
		
		\node at (1.5,-0.6) {\small no};

		
		\begin{scope}[xshift=4cm]
			\draw
			(0,0) -- (3,0)
			(0,-0.5) -- (0,0.5)
			(1.5,0) -- (1.5,0.5)
			(3,-0.5) -- (3,0.5)
			(2,0) -- (2,-0.5);
			\fill 
			(0,0) circle (0.05)
			(3,0) circle (0.05);
			
			\node at (0.15,0.2) {\small $\aaa$};
			\node at (1.35,0.2) {\small $\ccc$};
			\node at (0.8,0.2) {\small $b$};
			
			\node at (1.65,0.2) {\small $\ccc$};
			\node at (2.85,0.2) {\small $\aaa$};
			\node at (2.3,0.2) {\small $b$};
			
			\node at (0.15,-0.2) {\small $\bbb$};
			\node at (1.85,-0.2) {\small $\ccc$};
			\node at (2.15,-0.2) {\small $\bbb$};
			\node at (2.85,-0.2) {\small $\aaa$};
			\node at (1,-0.2) {$a$};
			\node at (2.5,-0.2) {$c$};
			
		\end{scope}
		
		
		\begin{scope}[xshift=8cm]
			\draw
			(0,0) -- (3,0)
			(0,-0.5) -- (0,0.5)
			(1.3,0) -- (1.3,0.5)
			(3,-0.5) -- (3,0.5)
			(1.7,0) -- (1.7,-0.5);
			\fill 
			(0,0) circle (0.05)
			(3,0) circle (0.05);
			
			\node at (0.15,0.2) {\small $\aaa$};
			\node at (1.15,0.2) {\small $\ccc$};
			\node at (0.15,-0.2) {\small $\bbb$};
			\node at (1.55,-0.2) {\small $\ccc$};
			\node at (1.85,-0.2) {\small $\ccc$};
			\node at (2.85,-0.2) {\small $\aaa$};
			\node at (0.65,0.2) {$b$};
			\node at (2.1,0.2) {$a$};
			\node at (0.9,-0.2) {$a$};
			\node at (2.3,-0.2) {$b$};
		\end{scope}

	\end{tikzpicture}
	\caption{Matching of $a+\cdots$ and $b+\cdots$ along an extended line}
	\label{ab_matching}
\end{figure}

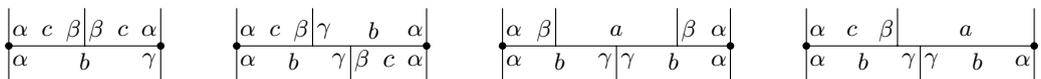
\begin{figure}[H]
	\centering
	\begin{tikzpicture}[>=latex,scale=1]
		
		\draw
		(0,0) -- (2,0)
		(0,-0.5) -- (0,0.5)
		(1,0) -- (1,0.5)
		(2,-0.5) -- (2,0.5);

		\fill 
		(0,0) circle (0.05)
		(2,0) circle (0.05);
		
		\node at (0.15,0.2) {\small $\aaa$};
		\node at (0.85,0.2) {\small $\bbb$};
		\node at (0.5,0.2) {\small $c$};
		
		\node at (1.15,0.2) {\small $\bbb$};
		\node at (1.85,0.2) {\small $\aaa$};
		\node at (1.5,0.2) {\small $c$};
		
		\node at (0.15,-0.2) {\small $\aaa$};
		\node at (1.85,-0.2) {\small $\ccc$};
		\node at (1,-0.2) {\small $b$};
		
		
		\begin{scope}[xshift=3cm]
			\draw
			(0,0) -- (2.5,0)
			(0,-0.5) -- (0,0.5)
			(1,0) -- (1,0.5)
			(2.5,-0.5) -- (2.5,0.5)
			(1.5,0) -- (1.5,-0.5);
			\fill 
			(0,0) circle (0.05)
			(2.5,0) circle (0.05);
			
			\node at (0.15,0.2) {\small $\aaa$};
			\node at (0.85,0.2) {\small $\bbb$};
			\node at (0.5,0.2) {\small $c$};
			
			\node at (1.15,0.2) {\small $\ccc$};
			\node at (2.35,0.2) {\small $\aaa$};
			\node at (1.8,0.2) {\small $b$};
			
			\node at (0.15,-0.2) {\small $\aaa$};
			\node at (1.35,-0.2) {\small $\ccc$};
			\node at (0.75,-0.2) {\small $b$};
			
			\node at (1.65,-0.2) {\small $\bbb$};
			\node at (2.35,-0.2) {\small $\aaa$};
			\node at (2,-0.2) {$c$};
			
		\end{scope}
		
		
		\begin{scope}[xshift=6.5cm]
			\draw
			(0,0) -- (3,0)
			(0,-0.5) -- (0,0.5)
			(0.7,0) -- (0.7,0.5)
			(2.3,0) -- (2.3,0.5)
			(1.5,0) --(1.5,-0.5)
			(3,-0.5) -- (3,0.5);
			\fill 
			(0,0) circle (0.05)
			(3,0) circle (0.05);
			
			\node at (0.15,0.2) {\small $\aaa$};
			\node at (0.55,0.2) {\small $\bbb$};

			\node at (1.5,0.2) {\small $a$};
			
			\node at (2.45,0.2) {\small $\bbb$};
			\node at (2.85,0.2) {\small $\aaa$};
			\node at (0.75,-0.2) {\small $b$};
			
			\node at (0.15,-0.2) {\small $\aaa$};
			\node at (1.35,-0.2) {\small $\ccc$};
			\node at (1.65,-0.2) {\small $\ccc$};
			\node at (2.85,-0.2) {\small $\aaa$};
			\node at (2.25,-0.2) {\small $b$};
		\end{scope}
		
		
		\begin{scope}[xshift=10.5cm]
			\draw
			(0,0) -- (3,0)
			(0,-0.5) -- (0,0.5)
			(1.2,0) -- (1.2,0.5)
			(3,-0.5) -- (3,0.5)
			(1.5,0) -- (1.5,-0.5);
			\fill 
			(0,0) circle (0.05)
			(3,0) circle (0.05);
			
			\node at (0.15,0.2) {\small $\aaa$};
			\node at (1.05,0.2) {\small $\bbb$};
			\node at (0.6,0.2) {\small $c$};
			
			\node at (2.1,0.2) {\small $a$};
			
			\node at (0.15,-0.2) {\small $\aaa$};
			\node at (1.35,-0.2) {\small $\ccc$};
			\node at (0.75,-0.2) {\small $b$};
			
			\node at (1.65,-0.2) {\small $\ccc$};
			\node at (2.85,-0.2) {\small $\aaa$};
			\node at (2.25,-0.2) {\small $b$};
			
		\end{scope}

	\end{tikzpicture}
	\caption{Matching of $b+\cdots$ and $c+\cdots$ along an extended line}
	\label{bc_matching}
\end{figure}
\begin{proof}
	In an extended edge, if there is only one stop vertex or no stop vertex, then it must be a great circle of length $2\pi$. If there are two stop vertices in an extended edge, then we have the following three possibilities.
	
Firstly, we start with two lengths $a$, $c$ at a vertex $\bullet$. Then the angles $\aaa$, $\bbb$ are at the end of $c$ as indicated in Figure \ref{ac_matching}. Consider the next length $x$ following $c$, i.e., $c+x+\cdots$. If $c+x\leq a$, then by $b+c>a$, we get $x=c$. If $2c=c+x<a$, then the $\aaa\bbb$-side $c$ implies $\aaa$ at a half vertex, a contradiction. Therefore, we have $x=c$ and $a=2c$ in the first picture of Figure \ref{ac_matching}.

If $c+x\geq a$, let $y$ be the next length following $a$, i,e., $a+y+\cdots$. If $x=b$, we have $y\geq c>x+c-a$ by $a>b>c$ in the second picture. This implies $y$ overhangs $\aaa$ on the right, a contradiction. Similarly, we have $x\neq c$. Therefore we have $x=a$. 

If $c+a\geq a+y$, we have $y=c$, get the third picture and $a+c=c+a$ along the extended edges. 

If $c+a< a+y$, then we have $y=b$ or $a$. let $z$ be the next length following $a_1$. If $y=b$, we get $c+a_1+z=a+b$ by no $\aaa$ at half vertex in the fourth picture. It implies $z=c=b-c$, and we get $a+b=c+a+c$ along the extended edges. 

We also need to consider the case $y=a$; the argument repeats by deleting a pair of $a$ from the two sides. We conclude that there are three types of extended edges: 
\begin{itemize}
	\item $(k+1)a=c+ka+c$.
	\item $c+ka=ka+c$.
	\item $c+ka+c=ka+b$.
\end{itemize}
    
    Secondly, we start with two lengths $a$, $b$ at a vertex $\bullet$. Then the angles $\aaa$, $\ccc$ at the end of $b$ as indicated in Figure \ref{ab_matching}. Consider the next length $x$ following $b$, i,e., $b+x+\cdots$. By $b+c>a$, we have $b+x>a$. Then let $y$ be the next length following $b$, i,e., $a+y+\cdots$.
   
   If $x=c$, we have $y\geq c>x+b-a$ in the first picture; this implies $y$ overhangs $\aaa$ on the right, a contradiction.
   
   If $x=b$, the extra length beyond $a$ is $b+x-a\leq c$. If $y=c$, we get the second picture and have $2b=a+c$ along the extended edges. If $y\neq c$, then the $\aaa\ccc$-side $b$ implies $\aaa$ at a half vertex, contradiction.
   
   If $x=a$, we have $y\neq c$ since there is no $\aaa$ at half vertex. If $y=b$, we get the third picture, and $b+a=a+b$ along the extended edges. 
   
   We also need to consider the case $y=a$; the argument repeats by deleting a pair of $a$ from the two sides. We conclude that there are two types of extended edges:
\begin{itemize}
	\item $b+ka+b=(k+1)a+c$.
	\item $b+ka=ka+b$.
\end{itemize}

Thirdly, we start with two lengths $b$, $c$ at a vertex $\bullet$. Then the angles $\aaa$, $\bbb$ at the end of $c$ as indicated in Figure \ref{bc_matching}. Consider the next length $x$ following $c$, i.e, $c+x+\cdots$, If $c+x\leq b$, then by $b+c>a$, we get $x=c$. If $2c=c+x<b$, then the $\aaa\bbb$-side implies $\aaa$ at a half vertex, a contradiction. Therefore we have either $x=c$ and $b=2c$ (first picture), or $x=c$ or $b$ with $c+x>b$	(second picture), or $x=a$ (third and fourth picture).

If $x=b$, we have $c+x>b$. Let $y$ be the next length following $b$, i.e., $b+y+\cdots$. The second picture's extra length beyond $b$ is $c+x-b\leq c$. Then, we get $y=c$ and $c+b=b+c$ along the extended edges by no $\aaa$ at the half vertex.

In the third and the fourth pictures, the extra length beyond $b$ is $c+x-b>c$; then, we get $y=b$ or $a$ by no $\aaa$ at the half vertex. If $y=b$ with $c+x<b+y=2b$, then the extra length $a$ is $2b-c-a\leq c$. Therefore, we get the third picture and $c+a+c=2b$ along the extended edges. If $c+x\geq b+y=2b$, we get the fourth picture and $c+a=2b$ by no $\aaa$ at the half vertex. 

We also need to consider the case $y=a$; the argument repeats by deleting a pair of $a$ from the two sides. We conclude that there are three types of extended edges: 

\begin{itemize}
	\item $c+ka+b=b+ka+c$.
	\item $c+(k+1)a+c=b+ka+b$.
	\item $c+(k+1)a=b+ka+b$.
\end{itemize}
\end{proof}
\end{lemma}

If $\aaa$ never appears at any half vertex, we also observe that any stop vertex must contain $\aaa$ from the proof of Lemma \ref{mis_M}. This fact induces the following lemma.

\begin{lemma}\label{extension edge}
	In a non-side-to-side $abc$-tiling, if $\aaa$ is not in any half vertex, then any half vertex is never a stop vertex for any extended edge. In other words, all corners at any half vertex are arranged side-to-side.
\end{lemma}
	
\begin{figure}[H]
		\centering
		\begin{tikzpicture}[>=latex,scale=0.8]
			\draw (4, 3)--(7, 3)
			(5.58, 2.96)--(5.56562, 5.56891)--(7.43548, 5.23401)--(5.55167, 4.32699)--(7.44943, 3.89441)--(5.58, 2.96)--(4, 5)--(5.56562, 5.56891);
			\node at (5.34235, 5.22006){\small $\ccc$};
			\node[inner sep=0.5, draw, shape=circle] at (5.11909, 4.29908) {\tiny 1};
			
			\node at  (6.77963, 5.16424){\small $\ccc$};
			\node[inner sep=0.5, draw, shape=circle] at (6, 5) {\tiny 3};
			
			
			\node at (6.79359, 3.85255){\small $\ccc$};
			\node[inner sep=0.5, draw, shape=circle] at (6.09588, 3.76883) {\tiny 2};
			
			\fill 
			(5.58, 2.96) circle (0.1);

		\end{tikzpicture}
		\caption{A possible half vertex when $\ccc$ is not in any half vertex.}
		\label{lemma 5 remark}
	\end{figure}
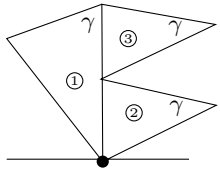

	\begin{remark}
	Similarly, if $\bbb$ is not in any half vertex, we also have all corners at any half vertex being arranged side-to-side. However, if $\ccc$ is not in any half vertex, we have no similar conclusion by the example in Figure \ref{lemma 5 remark}.
	\end{remark}
	
\section{Proof of the Main Theorem}
    In a non-side-to-side $abc$-tiling, a full vertex of degree 3,4,5 or a half vertex of degree 2 in Table \ref{vertex} must appear by Lemma \ref{345}. Then we can apply the irrational angle lemma to derive the AVC and to deduce the tiling. The main theorem will be proved by discussing all possible vertices in Table \ref{vertex}.  
    
   In the subsequent discussion, the tiles in a tiling will be labeled by circled numbers \raisebox{.2pt}{\textcircled{\raisebox{0.3pt} {\scriptsize $j$}}} in pictures, and denoted $T_j$. We add subscripts to indicate which tile certain corners and sides belong to. Therefore $\aaa_j$ and $a_j$ are the corner $\aaa$ and side $a$ in the tile $T_j$. 
   We use $\aaa\bbb^2\cdots$ to mean a vertex $\aaa^{m}\bbb^{n}\ccc^l$ with $m\geq1$ and $n\geq2$, and we use $R(\aaa\bbb^2\cdots)$	to mean the remainder. The remainder is either the remaining angle combination
   
   \[R(\aaa\bbb^2\cdots)=\aaa^{m-1}\bbb^{n-2}\ccc^l.
   \]
   or the remaining angle value
   \[R(\aaa\bbb^2\cdots)=2\pi-\aaa-2\bbb=(m-1)\aaa+(n-2)\bbb+l\ccc.
   \]
   We introduced these notations in our earlier work on spherical tilings \cite{wy1,wy2,awy}. 
    
\subsection{Tilings with a full vertex of degree 3}\label{vertex3}
    
    The AVC for $abc$-tilings with any full vertex of degree 3 are derived in Table \ref{V_3}, where AVC2 and AVC1 denote the collection of all other full vertices and all half vertices respectively.
    
    We will show the case $\aaa^3$ in full detail. The vertex $\aaa^3$ implies that $\aaa=\frac{2}{3}$ and $\bbb+\ccc=\frac{1}{3}+\frac{4}{f}>\frac{1}{3}$. By the irrational angle lemma, any full vertex $\aaa^{n_1}\bbb^{n_2}\ccc^{n_3}$ must satisfy $n_2 = n_3$. By $\bbb+\ccc>\frac{1}{3}$, we have $n_2=n_3\leq5$.  Then $\aaa^{n_1}\bbb^{n_2}\ccc^{n_3}=\aaa^3,\bbb^2\ccc^2,\aaa^2\bbb\ccc,\aaa\bbb^2\ccc^2,\bbb^3\ccc^3,\bbb^4\ccc^4,\aaa\bbb^3\ccc^3$ or $\bbb^5\ccc^5$.
    When $n_1,n_2,n_3$ are all even, $\aaa^{\frac{n_1}{2}}\bbb^{\frac{n_2}{2}}\ccc^{\frac{n_3}{2}}$ is a half vertex. These AVC are summarized in the first row of Table \ref{V_3}.
    All other AVC in Table \ref{V_3}, \ref{V_4}, \ref{V_5} are derived similarly.

\begin{table}[H]
	\centering
	\caption{The AVC for tilings with any full vertex of degree 3.}
	\begin{tabular}{|c|c|c|c|}
		\hline
		vertex &irrational angle lemma & AVC2 & AVC1 \\
		\hline
		\multirow{5}{*}{$\aaa^3$} &\multirow{5}{*}{$n_2=n_3$}  & $f=6$: $\bbb^2\ccc^2$&$f=6:\bbb\ccc$\\
		&&$f=12$: $\aaa^2\bbb\ccc$, $\aaa\bbb^2\ccc^2$, $\bbb^3\ccc^3$&\\
		&&$f=24$: $\bbb^4\ccc^4$&$f=24:\bbb^2\ccc^2$\\
		&&$f=36$: $\aaa\bbb^3\ccc^3$&\\
		&&$f=60$: $\bbb^5\ccc^5$&\\
		\hline
    	$\bbb^3$ &$n_1=n_3$ &$f=6$: $\aaa^2\ccc^2$&$f=6$: $\aaa\ccc$\\
    	\hline
    	$\aaa\bbb^2$&$2n_1=n_2+n_3$&$f=8$: $\aaa^2\ccc^4$&$f=8$: $\aaa\ccc^2$\\
		\hline
		\multirow{2}{*}{$\aaa^2\bbb$} & \multirow{2}{*}{$2n_2=n_1+n_3$} &$f=4k$: $\bbb^{\frac{f}{4}}\ccc^\frac{f}{2}$&$f=4k=8l:\bbb^\frac{f}{8}\ccc^\frac{f}{4}$\\
		&&$f=8k+4$: $\aaa\bbb^\frac{f+4}{8}\ccc^\frac{f}{4}$&\\
		\hline
		\multirow{2}{*}{$\aaa^2\ccc$} &\multirow{2}{*}{$2n_3=n_1+n_2$}  &$f=4k$: $\bbb^\frac{f}{2}\ccc^\frac{f}{4}$&$f=4k=8l:\bbb^\frac{f}{4}\ccc^\frac{f}{8}$\\
		&&$f=8k+4$: $\aaa\bbb^\frac{f}{4}\ccc^\frac{f+4}{8}$&\\
		\hline
		
	\end{tabular}

	\label{V_3}
\end{table}

\subsubsection*{Case $\aaa^3$}

From Table \ref{V_3}, this case has a half vertex only if $f=6,24$.
Case $f=6$ belongs to earth map tilings in the introduction. 

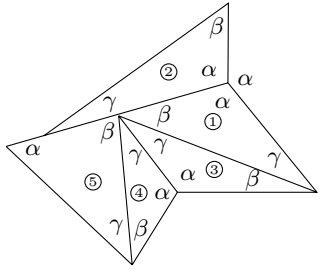
\begin{figure}[htp]
	\centering
	\begin{tikzpicture}[>=latex,scale=0.5]
    \draw (3.47826, 0.17031)--(4.66911, 2.09594)--(8.34302, 2.09594)--(6, 5)--(6.01199, 7.11273)--(1.14723, 3.59084)
    (0.15907, 3.31213)--(6, 5)
    (3.12354, 4.12292)--(3.47826, 0.17031)
    (3.12354, 4.12292)--(4.66911, 2.09594)
    (3.12354, 4.12292)--(8.34302, 2.09594)
    (3.47826, 0.17031)--(0.15907, 3.31213);
  
    \node at (5.85996, 4.47765) {\small $\aaa$};
    \node at (4.31439, 4.12292) {\small $\bbb$};
    \node at (7.22818, 2.95741){\small $\ccc$};
    \node[inner sep=0.5, draw, shape=circle] at (5.60659, 3.9709) {\tiny 1};

   \node at (5.50524, 5.31378) {\small $\aaa$};
   \node at (5.70794, 6.45396) {\small $\bbb$};
   \node at (2.8955, 4.42697){\small $\ccc$};
   \node[inner sep=0.5, draw, shape=circle] at (4.44108, 5.28844) {\tiny 2};

   \node at (4.94782, 2.55201){\small $\aaa$};
   \node at (6.65534, 2.40947) {\small $\bbb$};
   \node at (4.23838, 3.31213){\small $\ccc$};
   \node[inner sep=0.5, draw, shape=circle] at (5.64243, 2.70165) {\tiny 3};

   \node at (4.2789, 2.05884){\small $\aaa$};
   \node at (3.71401, 1.08489) {\small $\bbb$};
   \node at (3.57765, 3.03279){\small $\ccc$};
   \node[inner sep=0.5, draw, shape=circle] at (3.65557, 2.07832) {\tiny 4};
   
    \node at (0.87007, 3.2081){\small $\aaa$};
    \node at (2.83745, 3.65612) {\small $\bbb$};
    \node at (3.0712, 1.16281){\small $\ccc$};
    \node[inner sep=0.5, draw, shape=circle] at (2.44787, 2.37051) {\tiny 5};
    
 
    
     \node at (6.46055, 5.07809) {\small $\aaa$};

	\end{tikzpicture}
\caption{The possibilities for $\aaa^3$.}
\label{naaa-1}
\end{figure}

For $f=24$, we have $\text{AVC}=\{\aaa^3,\bbb^4\ccc^4;\,\bbb^2\ccc^2\}$. Since $\aaa$ is not in any half vertex, we have the AAD $\aaa^3=\thin^\ccc\aaa_1^\bbb\thin^\ccc\aaa_2^\bbb\thin\cdots$ in Figure \ref{naaa-1}, which  determines $T_1,T_2$. This gives a half vertex at $\bbb_1$. By Lemma \ref{extension edge}, we have $\bbb_1\cdots=\thin^\aaa\bbb_1^\ccc\thin^\bbb\ccc_3^\aaa\thin^\aaa\ccc_4^\bbb\thin^\ccc\bbb_5^\aaa\thin$, which determines $T_3$, $T_4$, $T_5$. Since $\aaa$ is not in any half vertex, we have $b\leq2c$. If $b<2c$, by $b+\text{all edge}>2c$, we get that $\aaa_5$ is a half vertex, a contradiction. Therefore $b=2c$, this gives a vertex $\aaa_5\ccc_2\cdots$, contradicting the AVC. 

\subsubsection*{Case $\bbb^3$}

 From Table \ref{V_3}, this case has a half vertex only if $f=6$. By $\bbb^3$, we have $\aaa+\ccc=1,\bbb=\frac23$. By Lemma \ref{Dawosn}, we have $\aaa+\bbb=\frac53-\ccc<1+\ccc$. This implies $\aaa<\frac23=\bbb$, contradicting $\aaa>\bbb$.

\subsubsection*{Case $\aaa\bbb^2$}
From Table \ref{V_3}, this case has a half vertex only if $f=8$.
By $\aaa\bbb^2$, we have $\aaa=1-2\ccc<1$, $\bbb=\frac{1}{2}+\ccc$. By $\aaa>\bbb$, we have $\ccc<\frac{1}{6}$. However, we have $\aaa+\bbb=\frac32-\ccc>1+\ccc$, contradicting Lemma \ref{Dawosn}.

\subsubsection*{Case $\aaa^2\bbb$}
From Table \ref{V_3}, this case has a half vertex only if $f=8l$.
For $f=8l(l\geq1)$,  $\text{AVC}=\{\aaa^2\bbb, \bbb^\frac{f}{4}\ccc^\frac{f}{2};\,\bbb^\frac{f}{8}\ccc^\frac{f}{4}\}$. Since $\aaa$ is not in any half vertex, we get the AAD $\aaa^2\bbb=\thin^\ccc\bbb^\aaa\thin^\bbb\aaa^\ccc\thin^\bbb\aaa^\ccc\thin$ or $\thin^\ccc\bbb^\aaa\thin^\bbb\aaa^\ccc\thin^\ccc\aaa^\bbb\thin$. In Figure \ref{2ab-1}, we have the AAD $\aaa^2\bbb=\thin^\ccc\bbb_1^\aaa\thin^\bbb\aaa_2^\ccc\thin^\bbb\aaa_3^\ccc\thin$, which determines $T_1,T_2,T_3$. There are half vertex at $\bbb_3$ and $\ccc_3$.  Since $\aaa$ is not in any half vertex, we derive $\aaa_4$.  By Lemma \ref{extension edge}, we have $\aaa_4\cdots=\aaa^3\cdots$, contradicting the AVC.

\begin{figure}[htp]
	\centering
	\begin{tikzpicture}[>=latex,scale=0.5]
		\draw (4, 0)--(4,-4)--(1,-7)--(8,-5)--(4,-4);
		\draw (8,-5)--(4,0);
		\draw (4,-1)--(1,-2)--(2,-6)--(4,-1);

		\node at (1.4, -2.25) {\small $\aaa$};
		\node at (1, -1.75) {\small $\aaa$};
		\node at (0.75, -2.4){\small $\aaa$};
		
		\node at (4.25, -1) {\small $\ccc$};
		\node at (3.75, -2.2) {\small $\bbb$};
		
		\node at (3.75, -3.8) {\small $\aaa$};
		\node at (4.25, -3.8) {\small $\aaa$};
		\node at (4, -4.4) {\small $\bbb$};
		
		\node at (2.6, -5.2) {\small $\ccc$};
		\node at (2.25, -6.3) {\small $\ccc$};
		
		\node at (7.1, -4.4) {\small $\bbb$};
		\node at (7.1, -5) {\small $\aaa$};
		
		\node[inner sep=0.5, draw, shape=circle] at (5.2, -3) {\tiny 2};
		\node[inner sep=0.5, draw, shape=circle] at (4.3, -5.5) {\tiny 1};
		\node[inner sep=0.5, draw, shape=circle] at (3.1, -4) {\tiny 3};
		\node[inner sep=0.5, draw, shape=circle] at (2.5, -3) {\tiny 4};
	
	\end{tikzpicture}
	\caption{The possibilities for $\aaa^2\bbb=\thin^\ccc\bbb^\aaa\thin^\bbb\aaa^\ccc\thin^\bbb\aaa^\ccc\thin$.}
	\label{2ab-1}
\end{figure}

Therefore, $\aaa^2\bbb=\thin^\ccc\bbb^\aaa\thin^\bbb\aaa^\ccc\thin^\ccc\aaa^\bbb\thin$.  Then $T_4$ in Figure \ref{2ab-2} is determined by $\aaa^2\bbb=\thin^\ccc\bbb^\aaa\thin^\bbb\aaa^\ccc\thin^\ccc\aaa^\bbb\thin$. If $\ccc_2\ccc_3\cdots$ is a full vertex, by Lemma \ref{mis_M}, we have $a=2c$. Then we can get $\aaa\ccc_2\cdots$, contradicting AVC. Therefore $\ccc_2\ccc_3\cdots$ is a half vertex. By the AVC and Lemma \ref{extension edge}, then $\ccc_2\ccc_3\cdots=\thin^\bbb\ccc^\aaa\thin^\aaa\ccc^\bbb\thin^\bbb\ccc^\aaa\thin^\aaa\ccc^\bbb\thin\cdots$,  $\thin^\bbb\ccc^\aaa\thin^\aaa\ccc^\bbb\thin^\ccc\bbb^\aaa\thin^\aaa\bbb^\ccc\thin\cdots$ or $\thin^\bbb\ccc^\aaa\thin^\aaa\ccc^\bbb\thin^\ccc\bbb^\aaa\thin$. 

If $\ccc_2\ccc_3\cdots=\thin^\bbb\ccc_2^\aaa\thin^\aaa\ccc_3^\bbb\thin^\bbb\ccc_6^\aaa\thin^\aaa\ccc_7^\bbb\thin\cdots$, then it determines $T_6$ and $T_7$ in the left of Figure \ref{2ab-2}. By Lemma \ref{extension edge}, we determine $T_8$. This gives a vertex $\aaa_6\aaa_7\aaa_8\cdots$, contradicting the AVC. Similarly, if $\ccc_2\ccc_3\cdots=\thin^\bbb\ccc^\aaa\thin^\aaa\ccc^\bbb\thin^\ccc\bbb^\aaa\thin^\aaa\bbb^\ccc\thin\cdots$ in the right of Figure \ref{2ab-2}, this also gives a vertex $\aaa_6\aaa_7\aaa_8\cdots$, contradicting the AVC.
\begin{figure}[H]
	\centering
	\begin{tikzpicture}[>=latex,scale=0.9]
		
		\draw (3.80054, 4.99003)--(7, 5)--(5.50975, 1.93537)--(5.50975, 3.93537)--(7, 5)
		(5.50975, 3.93537)--(3.80054, 4.99003)--(5.51,2.58)--(3.77967, 3.34277)--(3.80054, 4.99003)--(2.86214, 3.16325)--(3.77967, 3.34277)--(3.16134, 2.56486)--(5.51,2.58)
		(2.43537, 4.99322)--(5.2,4.99)--(3.41058, 6.13372)--(2.43537, 4.99322)
		
		(5.2,5) arc (135.56:57.92-182:2.02);
		
		\node at (6.49137, 4.41987) {\small $\aaa$};
		\node at (5.69352, 2.61473){\small $\ccc$};
		\node at (5.75336, 3.82148){\small $\bbb$};
		\node[inner sep=0.5, draw, shape=circle] at (6.0725, 3.65193) {\tiny 1};
		
		\node at (5.514, 4.18051){\small $\aaa$};
		\node at (4.59647, 4.82876){\small $\ccc$};
		\node at (6.35174, 4.79885){\small $\bbb$};
		\node[inner sep=0.5, draw, shape=circle] at (5.47411, 4.66919) {\tiny 2};
		
		\node at (5.31454, 3.83145){\small $\aaa$};
		\node at (4.49774, 4.34008){\small $\ccc$};
		\node at (5.3644, 3.19317){\small $\bbb$};
		\node[inner sep=0.5, draw, shape=circle] at (4.9954, 3.91124) {\tiny 3};
		
		\node at (7.22938, 5.07809){\small $\aaa$};
		\node at (5.85309, 2.05623){\small $\ccc$};
		\node at (5.69352, 5.1479){\small $\bbb$};
		\node[inner sep=0.5, draw, shape=circle] at (7.18949, 3.57215) {\tiny 4};
		
		\node[inner sep=0.5, draw, shape=circle] at (3.5, 5.5) {\tiny 5};
		
		\node at (4, 3.5){\small $\aaa$};
		\node at (3.97913, 4.37997){\small $\ccc$};
		\node at (4.96647, 3.08347){\small $\bbb$};
		\node[inner sep=0.5, draw, shape=circle] at (4.32819, 3.67188) {\tiny 6};
		
		\node at (3.60016, 3.48239){\small $\aaa$};
		\node at (3.63007, 4.2603){\small $\ccc$};
		\node at (3.24112, 3.45247){\small $\bbb$};
		\node[inner sep=0.5, draw, shape=circle] at (3.45056, 3.7018) {\tiny 7};
		
		\node at (3.85946, 3.1433){\small $\aaa$};
		\node at (3.55029, 2.73441){\small $\ccc$};
		\node at (4.79693, 2.75435){\small $\bbb$};
		\node[inner sep=0.5, draw, shape=circle] at (3.94921, 2.89398) {\tiny 8};	
		\begin{scope}[xshift=7cm]
			
			\draw (3.80054, 4.99003)--(7, 5)--(5.50975, 1.93537)--(5.50975, 3.93537)--(7, 5)
			(5.50975, 3.93537)--(3.80054, 4.99003)--(5.51,2.58)--(3.77967, 3.34277)--(3.80054, 4.99003)--(2.86214, 3.16325)--(3.77967, 3.34277)--(3.16134, 2.56486)--(5.51,2.58)
			(2.43537, 4.99322)--(5.2,4.99)--(3.41058, 6.13372)--(2.43537, 4.99322)
			
			(5.2,5) arc (135.56:57.92-182:2.02);
			
			\node at (6.49137, 4.41987) {\small $\aaa$};
			\node at (5.69352, 2.61473){\small $\ccc$};
			\node at (5.75336, 3.82148){\small $\bbb$};
			\node[inner sep=0.5, draw, shape=circle] at (6.0725, 3.65193) {\tiny 1};
			
			\node at (5.514, 4.18051){\small $\aaa$};
			\node at (4.59647, 4.82876){\small $\ccc$};
			\node at (6.35174, 4.79885){\small $\bbb$};
			\node[inner sep=0.5, draw, shape=circle] at (5.47411, 4.66919) {\tiny 2};
			
			\node at (5.31454, 3.83145){\small $\aaa$};
			\node at (4.49774, 4.34008){\small $\ccc$};
			\node at (5.3644, 3.19317){\small $\bbb$};
			\node[inner sep=0.5, draw, shape=circle] at (4.9954, 3.91124) {\tiny 3};
			
			\node at (7.22938, 5.07809){\small $\aaa$};
			\node at (5.85309, 2.05623){\small $\ccc$};
			\node at (5.69352, 5.1479){\small $\bbb$};
			\node[inner sep=0.5, draw, shape=circle] at (7.18949, 3.57215) {\tiny 4};
			
			\node[inner sep=0.5, draw, shape=circle] at (3.5, 5.5) {\tiny 5};
			
			\node at (4, 3.5){\small $\aaa$};
			\node at (4.96647, 3.08347){\small $\ccc$};
			\node at (3.97913, 4.37997){\small $\bbb$};
			\node[inner sep=0.5, draw, shape=circle] at (4.32819, 3.67188) {\tiny 6};
			
			\node at (3.60016, 3.48239){\small $\aaa$};
			\node at (3.24112, 3.45247){\small $\ccc$};
			\node at (3.63007, 4.2603){\small $\bbb$};
			\node[inner sep=0.5, draw, shape=circle] at (3.45056, 3.7018) {\tiny 7};
			
			\node at (3.85946, 3.1433){\small $\aaa$};
			\node at (4.79693, 2.75435){\small $\ccc$};
			\node at (3.55029, 2.73441){\small $\bbb$};
			\node[inner sep=0.5, draw, shape=circle] at (3.94921, 2.89398) {\tiny 8};		
		\end{scope}
		
	\end{tikzpicture}
	\caption{Two possibilities for $\thin^\ccc\bbb^\aaa\thin^\bbb\aaa^\ccc\thin^\ccc\aaa^\bbb\thin$.}
	\label{2ab-2}
\end{figure}

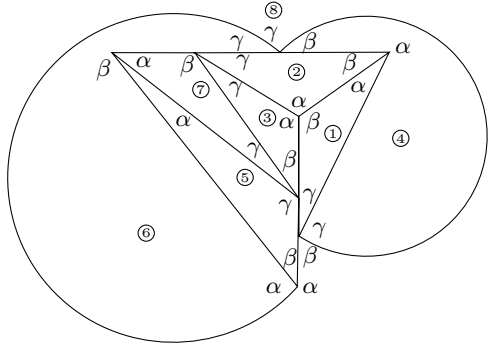
\begin{figure}[H]
	\centering
	\begin{tikzpicture}[>=latex,scale=0.8]				
		\draw (5.50975, 3.93537)--(5.50975, 1.93537)--(5.50975, 3.93537)--(3.80054, 4.99003)--(7, 5)--(5.50975, 3.93537)
		(5.51,2.58)--(3.80054, 4.99003)--(2.43537, 4.99322)-- (5.51,2.58)--(5.48508, 1.11876)--(2.43537, 4.99322)
		(7, 5)--(5.50975, 1.93537)
		(5.2,5) arc (135.56:57.92-182:1.99);
		\draw (5.2,5) arc (57.92-8:135.56+183:2.725);
		
		\node at (6.49137, 4.41987) {\small $\aaa$};
		\node at (5.69352, 2.61473){\small $\ccc$};
		\node at (5.75336, 3.82148){\small $\bbb$};
		\node[inner sep=0.5, draw, shape=circle] at (6.0725, 3.65193) {\tiny 1};
		
		\node at (5.514, 4.18051){\small $\aaa$};
		\node at (4.59647, 4.82876){\small $\ccc$};
		\node at (6.35174, 4.79885){\small $\bbb$};
		\node[inner sep=0.5, draw, shape=circle] at (5.47411, 4.66919) {\tiny 2};
		
		\node at (5.31454, 3.83145){\small $\aaa$};
		\node at (4.47779, 4.38995){\small $\ccc$};
		\node at (5.3644, 3.19317){\small $\bbb$};
		\node[inner sep=0.5, draw, shape=circle] at (4.9954, 3.91124) {\tiny 3};

		\node at (2.96187, 4.83874){\small $\aaa$};
		\node at (4.77699, 3.38266){\small $\ccc$};
		\node at (3.67994, 4.7789){\small $\bbb$};
		\node[inner sep=0.5, draw, shape=circle] at (3.89935, 4.37997) {\tiny 7};
		
		\node at (3.61013, 3.86137){\small $\aaa$};
		\node at (5.29559, 2.43521){\small $\ccc$};
		\node at (5.3654, 1.56755){\small $\bbb$};
		\node[inner sep=0.5, draw, shape=circle] at (4.62739, 2.9239) {\tiny 5};
		
		\node at (5.1, 1.1){\small $\aaa$};
		\node at (4.5, 5.15788){\small $\ccc$};
		\node at (2.3, 4.7){\small $\bbb$};
		\node[inner sep=0.5, draw, shape=circle] at (3, 2) {\tiny 6};
		
		\node at (7.22938, 5.07809){\small $\aaa$};
		\node at (5.85309, 2.05623){\small $\ccc$};
		\node at (5.69352, 5.1479){\small $\bbb$};
		\node[inner sep=0.5, draw, shape=circle] at (7.18949, 3.57215) {\tiny 4};

		
		\node at (5.05623, 5.2975){\small $\ccc$};
		\node at (5.7, 1.6){\small $\bbb$};
		\node at (5.7, 1.1){\small $\aaa$};
		\node[inner sep=0.5, draw, shape=circle] at (5.09613, 5.69643) {\tiny 8};
	\end{tikzpicture}
	\caption{Two possibilities of $\aaa^2\bbb$.} 
	\label{naab-1}
\end{figure}

Therefore, $\ccc_2\ccc_3\cdots=\thin^\bbb\ccc^\aaa\thin^\aaa\ccc^\bbb\thin^\ccc\bbb^\aaa\thin$. Then we have $f=8$ and $\text{AVC}=\{\aaa^2\bbb, \bbb^2\ccc^4;\,\bbb\ccc^2\}$. Therefore,  $T_4$ and $T_5$ are determined by the AVC in Figure \ref{naab-1}. By $\aaa^2\bbb=\thin^\ccc\bbb^\aaa\thin^\bbb\aaa^\ccc\thin^\ccc\aaa^\bbb\thin$, $T_6$ and $T_7$ are determined. This determines $T_8$. Then we get a unique tiling $T(4\aaa^2\bbb;\,4\bbb\ccc^2)$, whose 3D picture is shown in the first of Figure \ref{th-1}.  We have $\aaa=\frac{1}{2}+\ccc$, $\bbb=1-2\ccc$ and $\ccc\in(\frac{1}{4},\frac{1}{3})$.

\subsubsection*{Case $\aaa^2\ccc$}
We get a unique tiling $T(4\aaa^2\ccc;\,4\bbb^2\ccc)$ in the same way as $T(4\aaa^2\bbb;\,4\bbb\ccc^2)$ in Case $\aaa^2\bbb$, since $\bbb>\ccc$ has nothing to do with the deduction and we may switch $\bbb$ and $\ccc$. Then we have $\aaa=1-\frac{\ccc}{2}$, $\bbb=\frac{1}{2}-\frac{\ccc}{2}$ and $\ccc\in(\frac{1}{4},\frac{1}{3})$.

\vskip 9pt

Two tilings in Case $\aaa^2\bbb$ and Case $\aaa^2\ccc$ are essentially the same, and are combined into one tiling in the main theorem.

\vskip 9pt

After Section \ref{vertex3}, we can assume that no full vertex has degree $3$. 

\subsection{Tilings with a degree 4 full vertex or a degree 2 half vertex}\label{vertex4}
 The AVC for $abc$-tilings with any degree 4 full vertex or any degree 2 half vertex are derived in Table \ref{V_4}.

\begin{table}[H]
	\centering
	\caption{The AVC for tilings with a degree 4 full vertex or degree 2 half vertex.} 
	\begin{tabular}{|c|c|c|c|}
		\hline
		vertex
		&irrational angle lemma&AVC2&AVC1
		\\
		\hline
		$\aaa^3\bbb$
		&$3n_2=n_1+2n_3$&$f=12$: $\bbb^2\ccc^3$
		&none\\
		\hline
		$\aaa\bbb^3$
		&$3n_1=n_2+2n_3$&$f=12$: $\aaa^2\ccc^3$&none\\
		\hline
		$\bbb^3\ccc$
		&$3n_3=2n_1+n_2$&none&none\\
		\hline
		$\bbb\ccc^3$
		&$3n_2=2n_1+n_3$&none&none\\
		\hline
		$\aaa^3\ccc$
		&$3n_3=n_1+2n_2$&$f=12$: $\bbb^3\ccc^2$&none\\
		\hline	
		$\aaa\ccc^3$
		&$3n_1=2n_2+n_3$&none&none\\
		\hline
		\multirow{2}{*}{$\aaa^2\bbb\ccc$}
		&\multirow{2}{*}{$n_2=n_3$}&$f=4k$: $\bbb^\frac{f}{4}\ccc^\frac{f}{4}$&$f=4k=8l:\bbb^\frac{f}{8}\ccc^\frac{f}{8}$\\
		&&$f=8k+4$: $\aaa\bbb^\frac{f+4}{8}\ccc^\frac{f+4}{8}$&\\
		\hline
	    $\aaa\bbb^2\ccc$
		&$n_1=n_3$&$f=8$: $\bbb^4$, $\aaa^2\ccc^2$&$f=8$: $\bbb^2$, $\aaa\ccc$\\
		\hline
		\multirow{2}{*}{$\aaa\bbb\ccc^2$}
		&\multirow{2}{*}{$n_1=n_2$}&$f=5$: $\ccc^{10}$&$f=5$: $\ccc^5$\\
		&&$f=6$: $\ccc^8$&$f=6$: $\ccc^4$\\
		\hline
		\multirow{3}{*}{$\aaa^4$}
		&\multirow{3}{*}{$n_2=n_3$}&$f=8$: $\aaa^2\bbb\ccc$, $\bbb^2\ccc^2$&$f=$all: $\aaa^2$, $f=8$: $\bbb\ccc$\\
		&&$f=16$: $\aaa\bbb^2\ccc^2$&\\
		&&$f=24$: $\bbb^3\ccc^3$&\\
	   \hline
		\multirow{3}{*}{$\bbb^4$}
		&\multirow{3}{*}{$n_1=n_3$}&$f=8$: $\aaa^2\ccc^2$, $\aaa\bbb^2\ccc$&$f=$all: $\bbb^2$, $f=8$: $\aaa\ccc$\\
		&&$f=16$: $\aaa^2\bbb\ccc^2$&\\
		&&$f=24$: $\aaa^3\ccc^3$&\\
		\hline
		$\ccc^4$
		&$n_1=n_2$&none&$f=$all: $\ccc^2$\\
		\hline
		\multirow{2}{*}{$\aaa^2\bbb^2$}
		&\multirow{2}{*}{$n_1=n_2$}&$f=2k$: $\ccc^\frac{f}{2}$&$f=$all: $\aaa\bbb$\\
		&&$f=4k$: $\aaa\bbb\ccc^\frac{f}{4}$&$f=2k=4l$: $\ccc^\frac{f}{4}$\\
		\hline
		\multirow{2}{*}{$\aaa^2\ccc^2$}
		&\multirow{2}{*}{$n_1=n_3$}&$f=2k$: $\bbb^\frac{f}{2}$&$f=$all: $\aaa\ccc$\\
		&&$f=4k$: $\aaa\bbb^\frac{f}{4}\ccc$&$f=2k=4l$: $\bbb^\frac{f}{4}$\\
		\hline
		$\bbb^2\ccc^2$&$n_2=n_3$&none&$f=$all: $\bbb\ccc$\\
		\hline
	\end{tabular}
	
	\label{V_4}
\end{table}

The first six cases admit no half vertex. Case $\aaa\bbb\ccc^2$ and Case $\ccc^4$ cannot satisfy the balance lemma. So, we only need to discuss the other cases one by one. Clearly, Case $\aaa\bbb^2\ccc$, Case $\aaa^2\bbb^2$, $\aaa^2\ccc^2$ and  $\bbb^2\ccc^2$
in Table \ref{V_4} belong to earth map tilings in the introduction.

\subsubsection*{Case $\aaa^2\bbb\ccc$}
From Table \ref{V_4}, this case has a half vertex only if $f=8l$ and $\text{AVC}=\{\aaa^2\bbb\ccc,\bbb^\frac{f}{4}\ccc^\frac{f}{4};\,\bbb^\frac{f}{8}\ccc^\frac{f}{8}\}$.

For $l=1$, we have $\aaa=\frac12,\bbb+\ccc=1$. This implies $\aaa<\bbb$, a contradiction.
Therefore $l\ge2$. Then we get the AAD $\aaa^2\bbb\ccc=\thin^\aaa\ccc^\bbb\thin^\ccc\aaa^\bbb\thin^\ccc\bbb^\aaa\thin^\bbb\aaa^\ccc\thin$,
$\thin^\aaa\ccc^\bbb\thin^\bbb\aaa^\ccc\thin^\ccc\bbb^\aaa\thin^\bbb\aaa^\ccc\thin$, 
$\thin^\aaa\ccc^\bbb\thin^\ccc\aaa^\bbb\thin^\aaa\bbb^\ccc\thin^\ccc\aaa^\bbb\thin$, 
$\thin^\aaa\ccc^\bbb\thin^\ccc\aaa^\bbb\thin^\aaa\bbb^\ccc\thin^\bbb\aaa^\ccc\thin$ or
$\thin^\aaa\ccc^\bbb\thin^\ccc\bbb^\aaa\thin^\bbb\aaa^\ccc\thin\aaa\thin$.
If $l\ge3$, then we have the AAD $\bbb^{\aaa}\thin^{\aaa}\bbb\cdots$ at $\bbb^l\ccc^l$. This gives a vertex $\aaa^{\bbb}\thin^{\bbb}\aaa\cdots$, a contradiction. Therefore $l=2$.

In the left of Figure \ref{naabc-1}, we have the AAD $\aaa^2\bbb\ccc=\thin^\aaa\ccc^\bbb\thin^\ccc\aaa^\bbb\thin^\ccc\bbb^\aaa\thin^\bbb\aaa^\ccc\thin$, which determines $T_1$, $T_2$, $T_3$ and $T_4$. This gives half vertices at $\ccc_2$ and $\bbb_2$. By Lemma \ref{extension edge}, we can determines $\aaa_5$ and $\aaa_5\cdots=\aaa^3\cdots$, contradicting the AVC. 
Similarly, we have $\aaa^2\bbb\ccc\neq\thin^\aaa\ccc^\bbb\thin^\bbb\aaa^\ccc\thin^\ccc\bbb^\aaa\thin^\bbb\aaa^\ccc\thin$ and $\thin^\aaa\ccc^\bbb\thin^\ccc\aaa^\bbb\thin^\aaa\bbb^\ccc\thin^\ccc\aaa^\bbb\thin$.

In the right of Figure \ref{naabc-1}, we have the AAD  $\aaa^2\bbb\ccc=\thin^\aaa\ccc^\bbb\thin^\ccc\aaa^\bbb\thin^\aaa\bbb^\ccc\thin^\bbb\aaa^\ccc\thin$, which determines $T_1$, $T_2$, $T_3$ and $T_4$. By Lemma \ref{extension edge} and $\aaa\bbb\cdots=\aaa^2\bbb\ccc$, we can determines $T_5$ and $T_6$. This gives $\aaa^{\bbb}\thin^{\bbb}\aaa\cdots$, a contradiction. Therefore $\aaa^2\bbb\ccc=\thin^\aaa\ccc^\bbb\thin^\ccc\bbb^\aaa\thin^\bbb\aaa^\ccc\thin\aaa\thin$.
\begin{figure}[htp]
	\centering
	\begin{tikzpicture}[>=latex,scale=0.7]
		\draw (7, 7)--(10, 7)--(11.3, 10.22)--(10.76, 5.54)--(10, 7)
		(10.76, 5.54)--(7.4, 5.66)--(7, 7)--(10.76, 5.54)
		(10.68, 4.18)--(7.4, 5.66)--(11.48, 5.54)--(10.68, 4.18)
		(7.4, 5.66)--(3.6, 6.7)--(7, 7)
		(4.78,6.38) to[out=-90, in=180] (10.68, 4.18);
		
		\node at (10.63301, 4.49679)  {\small $\aaa$};
		\node at (11.0361, 5.26395){\small $\bbb$};
		\node at (8.34453, 5.43298){\small $\ccc$};
		\node[inner sep=0.5, draw, shape=circle] at (9.78783, 5.09491){\tiny 1};   
		
		\node at (7.66839, 5.95309)  {\small $\aaa$};
		\node at (7.53836, 6.48621){\small $\bbb$};
		\node at (9.47577, 5.81006){\small $\ccc$};
		\node[inner sep=0.5, draw, shape=circle] at (8.29252, 6.17414){\tiny 2};   
		
		\node at (6.86222, 6.74626)  {\small $\aaa$};
		\node at (7.07026, 6.13513){\small $\bbb$};
		\node at (5.19787, 6.59023){\small $\ccc$};
		\node[inner sep=0.5, draw, shape=circle] at (6.18607, 6.51221){\tiny 3};     
		
		\node at (7.2523, 5.39398)  {\small $\aaa$};
		\node at(5, 6){\small $\bbb$};
		\node at (9.85285, 4.34075){\small $\ccc$};
		\node[inner sep=0.5, draw, shape=circle] at(6.82321, 5.14692){\tiny 4};

		\node at (9.78665, 6.72144)  {\small $\aaa$};
		\node at(10.18808, 6.08456){\small $\bbb$};
		\node at (8.36792, 6.77176){\small $\ccc$};
		\node[inner sep=0.5, draw, shape=circle] at(9.35229, 6.62318){\tiny 5};

		\node at (10.26894, 7.00632) {\small $\aaa$};
		\node at (10.67202, 6.25216){\small $\bbb$};
		\node at (11, 9){\small $\ccc$};
		\node[inner sep=0.5, draw, shape=circle] at (10.68502, 7.63045) {\tiny 6};		

		\begin{scope}[xshift=-18cm]
			\draw (13, 4)--(17.30153, 3.98581)--(18.89881, 5.54594)--(14.40415, 5.65738)--(17.30153, 3.98581)
			(14.40415, 5.65738)--(13, 9)--(13, 4)
			(13.36,8.15)--(17.78,5.57)--(16.54004, 8.12759)--(13.36,8.15)
			(14.40415, 5.65738)--(13, 4);
			
			\node at (17.26477, 4.46219) {\small $\aaa$};
			\node at (18.19895, 5.26526){\small $\bbb$};
			\node at (15.5439, 5.34721){\small $\ccc$};
			\node[inner sep=0.5, draw, shape=circle] at (16.74031, 5.01943){\tiny 1}; 
			
			\node at (14.64249, 5.93722) {\small $\aaa$};
			\node at(13.98692, 7.44503){\small $\bbb$};
			\node at (16.67476, 5.90444){\small $\ccc$};
			\node[inner sep=0.5, draw, shape=circle] at (15.2325, 6.49446){\tiny 2};   
			
			\node at (13.33135, 4.77359) {\small $\aaa$};
			\node at(14.06887, 5.675){\small $\bbb$};
			\node at (13.2658, 7.67448){\small $\ccc$};
			\node[inner sep=0.5, draw, shape=circle] at (13.46247, 6.18306){\tiny 3};       
			
			\node at (14.4786, 5.31443) {\small $\aaa$};
			\node at(13.59358, 4.28191){\small $\bbb$};
			\node at (16.21586, 4.26552){\small $\ccc$};
			\node[inner sep=0.5, draw, shape=circle] at (14.72444, 4.7572){\tiny 4};
			
			\node at (16.37975, 7.77282) {\small $\aaa$};
			\node at (16.9, 8) {\small $\aaa$};
			\node at (16.2, 8.3) {\small $\aaa$};
			
			\node[inner sep=0.5, draw, shape=circle] at (16.1503, 7.3467){\tiny 5};

		\end{scope}
		
	\end{tikzpicture}
	\caption{Two possibilities for $\aaa^2\bbb\ccc$.}
	\label{naabc-1}
\end{figure}

If $\aaa^2\bbb\ccc=\thin^\aaa\ccc_1^\bbb\thin^\ccc\bbb_2^\aaa\thin^\bbb\aaa_3^\ccc\thin^\bbb\aaa_4^\ccc\thin$, then it determines $T_1$, $T_2$, $T_3$ and $T_4$ in Figure \ref{naabc-2}. This gives a half vertex at $\bbb_4$. By Lemma \ref{extension edge}, we have $\bbb_4\cdots=\thin\bbb_4^{\ccc}\thin^{\bbb}\ccc_5^{\aaa}\thin^{\aaa}\ccc_6^{\bbb}\thin^{\ccc}\bbb_7\thin$. This determine $T_5,T_6,T_7$ and $b=2c$. Then $\thin^{\bbb}\ccc_3^{\aaa}\thin^{\bbb}\aaa_7^{\bbb}\thin\cdots=\thin^{\bbb}\ccc_3^{\aaa}\thin^{\bbb}\aaa_7^{\bbb}\thin^{\ccc}\aaa_8^{\bbb}\thin^{\aaa}\bbb_9^{\ccc}\thin$ determine $T_8,T_9$;  $\thin^{\bbb}\aaa_1^{\ccc}\thin^{\aaa}\ccc_4^{\bbb}\thin^{\ccc}\bbb_5^{\aaa}\thin\cdots=\thin^{\bbb}\aaa_1^{\ccc}\thin^{\aaa}\ccc_4^{\bbb}\thin^{\ccc}\bbb_5^{\aaa}\thin^{\bbb}\aaa_{10}^{\ccc}\thin$ determines $T_{10}$; This gives a half vertex at $\bbb_1$. By Lemma \ref{extension edge}, we have $\thin^{\aaa}\bbb_1^{\ccc}\thin^{\bbb}\ccc_2^{\aaa}\thin\cdots=\thin^{\aaa}\bbb_1^{\ccc}\thin^{\bbb}\ccc_2^{\aaa}\thin^{\aaa}\ccc_{12}^{\bbb}\thin^{\ccc}\bbb_{13}^{\aaa}\thin$ determine $T_{12}, T_{13}$. Similarly, we determine $T_{14}$, $T_{15}$ and $T_{16}$. Then we get a unique tiling $T(8\aaa^2\bbb\ccc;4\bbb^2\ccc^2)$, the 3D picture for tiling is shown in the third of Figure \ref{th-1}. If $\aaa^2\bbb\ccc=\thin^\aaa\ccc_1^\bbb\thin^\ccc\bbb_2^\aaa\thin^\bbb\aaa_3^\ccc\thin^\ccc\aaa_4^\bbb\thin$, we also obtain this tiling.

\begin{figure}[H]
	\centering
	\begin{tikzpicture}[>=latex,scale=0.8]
		\draw (7.41491, 4.97133)--(9, 5)--(8.99301, 3.26652)--(7.41491, 4.97133)--(7.43795, 3.255)--(8.99301, 3.26652)--(5.79074, 3.24348)--(7.41491, 4.97133)--(5.80226, 4.97133)--(5.79074, 3.24348)--(7.41491, 4.97133)
		(5.80226, 4.97133)--(4.59647, 3.23235)--(4.55266, 5.7403)--(3.42463, 5.75126)--(4.59647, 3.23235)--(5.79074, 3.24348)
		(5.80226, 4.97133)--(4.55266, 5.7403)--(9, 5)
		(4.59647, 3.23235)--(3.45748, 3.2433)--(3.42463, 5.75126)--(2.35135, 3.27615)--(3.45748, 3.2433)
		(2.35135, 3.27615)--(2.32945, 5.76221)--(3.42463, 5.75126)
		(2.32945, 5.76221) arc (120.56:57.92-11.5:5.575)
		(2.32945, 5.76221) arc (164.56:57.92-83.8:3.575)
		(2.33945, 3.28235) arc (-140:-40.5:4.36)
		(3.45748, 3.2433) to[out=-90, in=-90](7.43795, 3.255);

		\node at (4.76772, 5.42115) {\small $\aaa$};
		\node at(5.39863, 4.80042){\small $\bbb$};
		\node at (4.7779, 3.78284){\small $\ccc$};
		\node[inner sep=0.5, draw, shape=circle] at (5.02212, 4.68849){\tiny 2};
		
		\node at (5.61232, 3.40633) {\small $\aaa$};
		\node at(4.96106, 3.40633){\small $\bbb$};
		\node at (5.67337, 4.48497){\small $\ccc$};
		\node[inner sep=0.5, draw, shape=circle] at (5.46986, 3.9253){\tiny 1};
		
		\node at (5.999, 4.73937) {\small $\aaa$};
		\node at (6.85378, 4.70884){\small $\bbb$};
		\node at (5.96847, 3.65055){\small $\ccc$};
		\node[inner sep=0.5, draw, shape=circle] at (6.32463, 4.42392){\tiny 4};
		
		\node at (5.89724, 5.13623) {\small $\aaa$};
		\node at (5.51056, 5.38045){\small $\bbb$};
		\node at (7.29134, 5.13623){\small $\ccc$};
		\node[inner sep=0.5, draw, shape=circle] at (6.40604, 5.22869){\tiny 3};
		
		\node at (7.26081, 3.39615) {\small $\aaa$};
		\node at (6.22287, 3.43686){\small $\bbb$};
		\node at (7.24046, 4.46462){\small $\ccc$};
		\node[inner sep=0.5, draw, shape=circle] at (7, 4){\tiny 5};
		
		\node at (7.65767, 3.43686) {\small $\aaa$};
		\node at (8.60403, 3.44703){\small $\bbb$};
		\node at (7.63732, 4.48497){\small $\ccc$};
		\node[inner sep=0.5, draw, shape=circle] at (7.87136, 4.01688){\tiny 6};
		
		\node at (5.79549, 3.00947) {\small $\aaa$};
		\node at (7.1794, 3.05017){\small $\bbb$};
		\node at (3.78066, 3.06035){\small $\ccc$};
		\node[inner sep=0.5, draw, shape=circle] at (5.56144, 2.71437)  {\tiny 10};
		
		\node at (4.36069, 5.54326) {\small $\aaa$};
		\node at (3.83154, 5.52291) {\small $\bbb$};
		\node at (4.43192, 4.07794) {\small $\ccc$};
		\node[inner sep=0.5, draw, shape=circle] at (4.24875, 4.91236)  {\tiny 12};
		
		\node at (2.90026, 5.88) {\small $\aaa$};
		\node at (8.3, 5.3) {\small $\bbb$};
		\node at (4.65579, 6.01135) {\small $\ccc$};
		\node[inner sep=0.5, draw, shape=circle] at (5.55968, 6.13698)  {\tiny 9};
		\node[inner sep=0.5, draw, shape=circle] at (5.55968, 6.8698)  {\tiny 8};
		\node at (2.7, 6.2) {\small $\bbb$};
		\node at (2.1, 5.7) {\small $\ccc$};
		
		\node at (3.64838, 3.42668) {\small $\aaa$};
		\node at (4.2691, 3.4165) {\small $\bbb$};
		\node at (3.5975, 4.97341) {\small $\ccc$};
		\node[inner sep=0.5, draw, shape=circle] at (3.80101, 4.28145)  {\tiny 13};
		
		\node at (3.25152, 3.46738) {\small $\aaa$};
		\node at (2.62061, 3.43686) {\small $\bbb$};
		\node at (3.27187, 4.89201) {\small $\ccc$};
		\node[inner sep=0.5, draw, shape=circle] at (3.12941, 4.29163)  {\tiny 14};
		
		\node at (2.1, 3.3) {\small $\aaa$};
		\node at (2.75, 3.1) {\small $\aaa$};
		\node at (2.51885, 5.54326) {\small $\aaa$};
		\node at (3.10905, 5.56361) {\small $\bbb$};
		\node at (2.4985, 4.09829) {\small $\ccc$};
		\node[inner sep=0.5, draw, shape=circle] at (2.65114, 4.83095)  {\tiny 15};
		\node[inner sep=0.5, draw, shape=circle] at (2, 4.83095)  {\tiny 16};
		
		\node at (8.8279, 4.78007) {\small $\aaa$};
		\node at (7.9833, 4.7699) {\small $\bbb$};
		\node at (8.83807, 3.7116){\small $\ccc$};
		\node[inner sep=0.5, draw, shape=circle] at (8.49209, 4.40356) {\tiny 7};
		
		\node at (9.18, 4.96324) {\small $\aaa$};
		\node at (9.15, 4.0324) {\small $\ccc$};
		\node at (9.1, 3.1) {\small $\bbb$};
		\node at (3.28204, 3.04) {\small $\bbb$};
		\node at (7.64749, 3.01965){\small $\ccc$};
		\node[inner sep=0.5, draw, shape=circle] at (7.8, 2.63296) {\tiny 11};

	\end{tikzpicture}\hspace{11pt}

	\caption{The tiling for $\aaa^2\bbb\ccc=\thin^\aaa\ccc^\bbb\thin^\ccc\bbb^\aaa\thin^\bbb\aaa^\ccc\thin^\bbb\aaa^\ccc\thin$.}
	\label{naabc-2}
\end{figure}

By calculation, we get 
\[\aaa=\frac{3}{4},\bbb=\frac{\arctan\sqrt{2}}{\pi}\approx0.3041,\ccc=\frac{1}{2}-\frac{\arctan\sqrt{2}}{\pi}\approx0.1959,\]
\[a=\frac{2\pi}{3},b=\frac{\pi}{2},c=\frac{\pi}{4}.\]
 
 \subsubsection*{Case $\aaa^4$}
 From Table \ref{V_4}, there exists a half vertex only if $f=8,16,24$.
 
 For $f=8$, we have $\aaa=\frac12,\bbb+\ccc=1$. This implies $\bbb>\frac12=\aaa$, contradiction.

\begin{figure}[H]
 	\centering
 	\includegraphics[width=0.8\textwidth]{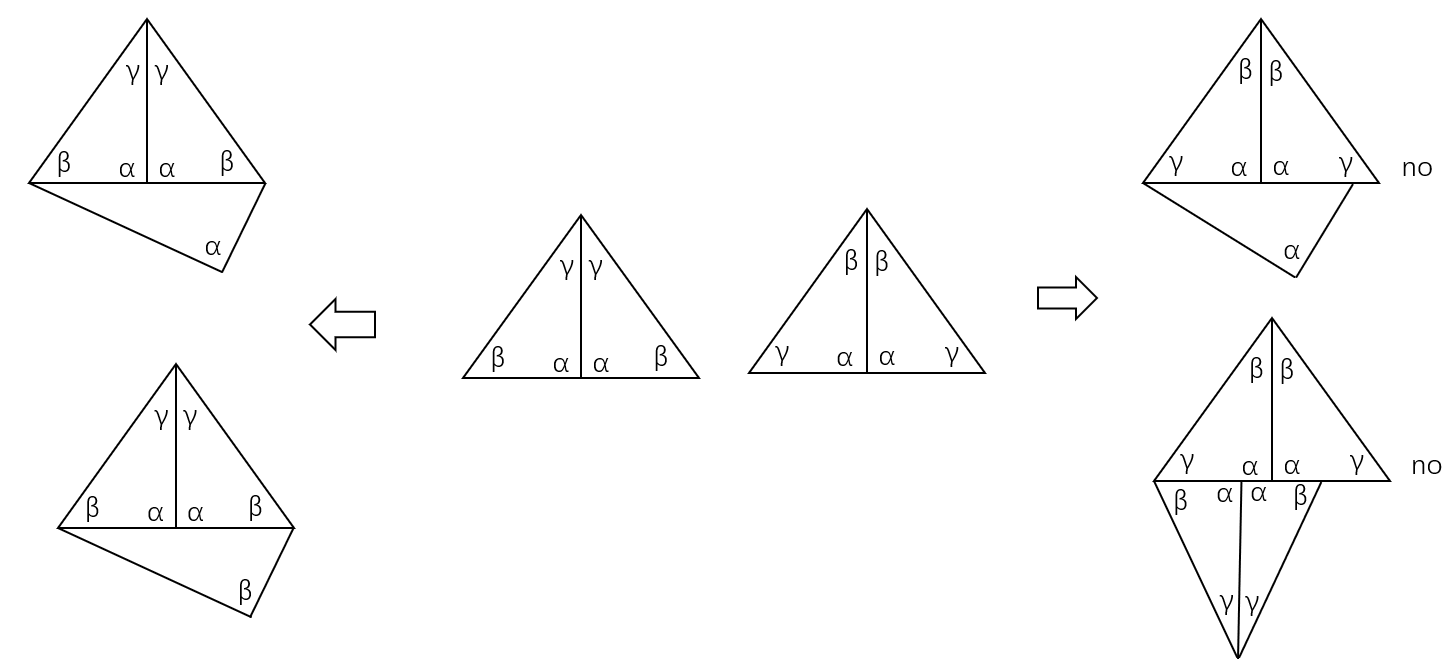}
 	\caption{The possibilities for a half vertex $\aaa\thin\aaa$.}
 	\label{naaaa-0}
 \end{figure}

For $f=16,24$, $\bbb$ and $\ccc$ are not in any half vertex by the AVC. This implies that all half vertices are $\aaa^2$. By AAD, we have the half vertex $\aaa^\ccc\thin\aaa=\aaa^\ccc\thin^\ccc\aaa$ or $\aaa^\bbb\thin\aaa=\aaa^\bbb\thin^\bbb\aaa$. Then we conclude the two AADs of $\aaa\thin\aaa$ in the Figure \ref{naaaa-0}. For the left of Figure \ref{naaaa-0}, we get $b=2c$ or $a=2c$. For the right of Figure \ref{naaaa-0}, if $a$ overhangs $\aaa^2$, then this gives a half vertex at $\bbb$ or $\ccc$ by $2b>a$, a contradiction. Similarly, $c$ can not overhang $\aaa^2$. Therefore, the half vertex $\aaa\thin\aaa$ must be $\aaa^\ccc\thin^\ccc\aaa$ and  $b=2c$ or $a=2c$.

For $f=16$, we have $\text{AVC}=\{\aaa^4,\aaa\bbb^2\ccc^2;\,\aaa^2\}$. If $a=2c$, in the adjacent tiles of $\aaa\bbb^2\ccc^2$, there is a half vertex at $\bbb$ or $\ccc$.
If $b=2c$, then we get $\aaa=\frac{1}{2},\bbb=\frac{1}{2},\ccc=\frac{1}{4}$ by calculation. Dawson has discussed in \cite{dawson2001}.

For $f=24$, we have $\text{AVC}=\{\aaa^4,\bbb^3\ccc^3;\,\aaa^2\}$. If $a=2c$,  in the adjacent tiles of $\aaa\bbb^2\ccc^2$, there is a half vertex at $\bbb$ or $\ccc$. If $b=2c$, then there is a vertex $\aaa\bbb\cdots$, contradicting the AVC.

\subsubsection*{Case $\bbb^4$}
 From Table \ref{V_4}, there exists a half vertex only if $f=8,16,24$.

 For $f=8$, it admits earth map tilings.
  
 For $f=16,24$, $\aaa$ and $\ccc$ are not in any half vertex by the AVC. This implies that all half vertices are $\bbb^2$. By AAD, we have any half vertex is either $\bbb^\ccc\thin^\ccc\bbb$ or $\bbb^\aaa\thin^\aaa\bbb$. Similar to Case $\aaa^4$, we derive that any half vertex must be $\bbb^\ccc\thin^\ccc\bbb$ and $a=2c$ or $b=2c$.

 For $f=16$, then we have $\aaa+\ccc=\frac{3}{4}, \bbb=\frac{1}{2}$. If $a=2c$, then we get this triangle is a isosceles triangle ($\frac{1}{2},\frac{1}{2},\frac{1}{4}$) by calculation. Dawson has discussed in \cite{dawson2001}.  If $b=2c$,
 then we get $a=\frac{1}{2}, c=\frac{3}{4}, b=\frac{3}{2}$, contradicting $a>c$.

 For $f=24$, we have $\text{AVC}=\{\aaa^3\ccc^3,\bbb^4;\,\bbb^2\}$. If $a=2c$, this gives a vertex $\aaa\bbb\cdots$ by AAD, contradicting the AVC. If $b=2c$, similar to Case $\aaa^4$, we get a similar contradiction in a full vertex $\aaa^3\ccc^3$. 

\vskip 9pt
 
After Section \ref{vertex4}, we can assume that no full vertex has degree $3$ or $4$, and no half vertex has degree $2$.

 \subsection{Tilings with a full vertex of degree 5} 
  The AVC for $abc$-tilings with any full vertex of degree 5 are derived in Table \ref{V_5}. 
  
  In Table \ref{V_5},  the first $16$ cases up to Case $\aaa^3\bbb\ccc$ and the last case admit no half vertex. Case $\aaa\bbb\ccc^3$ and Case $\aaa\bbb^3\ccc$ cannot satisfy the balance lemma.

\begin{table}[H]
	\caption{The AVC for tilings with any full vertex of degree 5.} 
	\centering
	\begin{tabular}{|c|c|c|c|}
		\hline
		vertex
		&irrational angle lemma &AVC2&AVC1
		\\
		\hline
		$\aaa^5$
		&$n_2=n_3$&$f=60$: $\bbb^3\ccc^3$&none\\
		\hline
		\multirow{3}{*}{$\bbb^5$}
		&\multirow{3}{*}{$n_1=n_3$}&$f=10$: $\aaa^2\ccc^2$&\multirow{3}{*}{none}\\
		&&$f=20$: $\aaa\bbb^3\ccc$, $\aaa^2\bbb\ccc^2$&\\
		&&$f=60$: $\aaa^3\ccc^3$&\\
		\hline
		$\ccc^5$
		&$n_1=n_2$&none&none\\
		\hline
		$\aaa^4\bbb$
		&$4n_2=n_1+3n_3$&none&none\\
		\hline
		$\aaa^4\ccc$
		&$4n_3=n_1+3n_2$&none&none\\
		\hline
		$\aaa^3\bbb^2$
		&$3n_2=2n_1+n_3$&none&none\\
		\hline
		$\aaa^2\bbb^3$
		&$3n_1=2n_2+n_3$&none&none\\
		\hline
		$\aaa^3\ccc^2$
		&$3n_3=2n_1+n_2$&none&none\\
		\hline
		$\aaa^2\ccc^3$
		&$3n_1=n_2+2n_3$&none&none\\
		\hline
		$\aaa\bbb^4$
		&$4n_1=n_2+3n_3$&none&none\\
		\hline
		$\aaa\ccc^4$
		&$4n_1=3n_2+n_3$&none&none\\
		\hline
		$\bbb^4\ccc$
		&$4n_3=3n_1+n_2$&none&none\\
		\hline
		$\bbb\ccc^4$
		&$4n_2=3n_1+n_3$&none&none\\
		\hline
		$\bbb^3\ccc^2$
		&$3n_3=n_1+2n_2$&none&none\\
		\hline
		$\bbb^2\ccc^3$
		&$3n_2=n_1+2n_3$&none&none\\
		\hline
		$\aaa^3\bbb\ccc$
		&$n_2=n_3$&$f=36$: $\bbb^3\ccc^3$&none\\
		\hline
		\multirow{2}{*}{$\aaa^2\bbb^2\ccc$}
		&\multirow{2}{*}{$n_1=n_2$}&$f=4k$: $\ccc^\frac{f}{4}$&$f=4k=8l:\ccc^\frac{f}{8}$\\
		&&$f=8k+4$: $\aaa\bbb\ccc^\frac{f+4}{8}$&\\
		\hline
		\multirow{2}{*}{$\aaa^2\bbb\ccc^2$}
		&\multirow{2}{*}{$n_1=n_3$}&$f=4k$: $\bbb^\frac{f}{4}$&$f=4k=8l:\bbb^\frac{f}{8}$\\
		&&$f=8k+4$: $\aaa\bbb^\frac{f+4}{8}\ccc$&\\
		\hline
		\multirow{2}{*}{$\aaa\bbb\ccc^3$}
		&\multirow{2}{*}{$n_1=n_2$}&$f=8$: $\ccc^8$&$f=8$: $\ccc^4$\\
		&&$f=12$: $\ccc^6$&$f=12$: $\ccc^3$\\
		\hline
		\multirow{4}{*}{$\aaa\bbb^3\ccc$}
		&\multirow{4}{*}{$n_1=n_3$}&$f=8$: $\bbb^8$&$f=8$: $\bbb^4$\\
		&&$f=12$: $\bbb^6$&$f=12$: $\bbb^3$\\
		&&$f=20$: $\aaa^2\bbb\ccc^2,\bbb^5$&\\
		&&$f=36$: $\aaa^3\ccc^3$&\\
        \hline
		\multirow{2}{*}{$\aaa\bbb^2\ccc^2$}
		&\multirow{2}{*}{$n_2=n_3$}&$f=8$: $\bbb^4\ccc^4$&\multirow{2}{*}{none}\\
		&&$f=12$: $\bbb^3\ccc^3$&\\
		\hline
	\end{tabular}
	
	\label{V_5}
\end{table}

\subsubsection*{Case $\aaa^2\bbb^2\ccc$}
 From Table \ref{V_5}, this case has a half vertex only if $f=8l$. For $f=8l(l\geq3)$, we have $\text{AVC}=\{\aaa^2\bbb^2\ccc,\ccc^\frac{f}{4};\,\ccc^\frac{f}{8}\}$. This implies $\aaa$ and $\bbb$ are not in any half vertex. Then we have the AAD $\aaa^2\bbb^2\ccc=\thin^\bbb\ccc^\aaa\thin^{\ccc}\aaa^{\bbb}\thin^{\aaa}\bbb^{\ccc}\thin^{\ccc}\bbb^{\aaa}\thin^{\bbb}\aaa^{\ccc}\thin,\thin^\bbb\ccc^\aaa\thin^\ccc\aaa^\bbb\thin^\bbb\aaa^\ccc\thin^\ccc\bbb^\aaa\thin^\aaa\bbb^\ccc\thin$ or $\thin^\bbb\ccc^\aaa\thin^\ccc\aaa^\bbb\thin^\aaa\bbb^\ccc\thin^\ccc\aaa^\bbb\thin^\aaa\bbb^\ccc\thin$.

If $\aaa^2\bbb^2\ccc=\thin^\bbb\ccc_1^\aaa\thin^{\ccc}\aaa_2^{\bbb}\thin^{\aaa}\bbb_3^{\ccc}\thin^{\ccc}\bbb_4^{\aaa}\thin^{\bbb}\aaa_5^{\ccc}\thin$, this determine $T_1$, $T_2$, $T_3$, $T_4$ and $T_5$ in the first of Figure \ref{naabbc-1}. We get $\bbb_1\cdots$ is a half vertex by Lemma \ref{mis_M}, a contradiction.

If $\aaa^2\bbb^2\ccc=\thin^\bbb\ccc_1^\aaa\thin^\ccc\aaa_2^\bbb\thin^\bbb\aaa_3^\ccc\thin^\ccc\bbb_4^\aaa\thin^\aaa\bbb_5^\ccc\thin$, this determine $T_1$, $T_2$, $T_3$, $T_4$ and $T_5$ in the second of Figure \ref{naabbc-1}. This gives a half vertex at $\ccc_4$.  By Lemma \ref{extension edge}, 
$\ccc_4^{\bbb}\thin\cdots=\ccc_4^{\bbb}\thin^\bbb\ccc_5\thin\cdots$ determines $T_6$. This gives a vertex $\bbb^3\cdots$, contradicting the AVC.  Therefore $\aaa^2\bbb^2\ccc=\thin^\bbb\ccc^\aaa\thin^\ccc\aaa^\bbb\thin^\aaa\bbb^\ccc\thin^\ccc\aaa^\bbb\thin^\aaa\bbb^\ccc\thin$

We have the AAD $\aaa^2\bbb^2\ccc=\thin^\bbb\ccc_1^\aaa\thin^\ccc\aaa_2^\bbb\thin^\aaa\bbb_3^\ccc\thin^\ccc\aaa_4^\bbb\thin^\aaa\bbb_5^\ccc\thin$, this determine $T_1$, $T_2$, $T_3$, $T_4$ and $T_5$ in the third of Figure \ref{naabbc-1}. This gives a half vertex at $\ccc_5$. By Lemma \ref{extension edge}, $\ccc_5^{\bbb}\thin\cdots=\ccc_5^{\bbb}\thin^{\bbb}\ccc_6\cdots$ determines $T_6$. This gives a vertex $\aaa_3\thin\bbb_5\thin\bbb_6\cdots=\aaa^2\bbb^2\ccc$, a contradiction.

\begin{figure}[htp]
	\centering
	\begin{tikzpicture}[>=latex,scale=0.5]
		\draw (6.99539, 5.15148)--(10.02255, 5.13546)--(7, 7)
		(6.99539, 5.133)--(7,7)
		(6.99539, 5.133)--(3.54613, 7.00203)
		(10.02255, 5.13546)--(8.34951, 9.38681)--(7, 7)--(5.37277, 10.08025)--(3.54613, 7.00203)
		(8.34951, 9.38681)--(5.71104, 9.38681)
		(7, 7)--(3.54613, 7.00203);	
		
		\node at (6.70709, 6.73714){\small $\aaa$};
		\node at (6.65904, 5.80817){\small $\bbb$};
		\node at (4.76906, 6.75315){\small $\ccc$};
		\node[inner sep=0.5, draw, shape=circle] at (5.92227, 6.35273){\tiny 2};
		
		\node at (7.29971, 5.48783){\small $\aaa$};
		\node at (7.23564, 6.4){\small $\bbb$};
		\node at (9.06155, 5.48783){\small $\ccc$};
		\node[inner sep=0.5, draw, shape=circle] at (7.78021, 6.01638){\tiny 3};
		
		\node at (8.30876, 8.62711){\small $\aaa$};
		\node at (7.55598, 7.18561){\small $\bbb$};
		\node at (9.38188, 5.90427){\small $\ccc$};
		\node[inner sep=0.5, draw, shape=circle] at (8.38885, 7.28171){\tiny 4};
		
		\node at (7.01141, 7.58602){\small $\aaa$};
		\node at (7.73216, 9.07558){\small $\bbb$};
		\node at (6.21057, 9.10762){\small $\ccc$};
		\node[inner sep=0.5, draw, shape=circle] at (7.01141, 8.35483){\tiny 5};
		
		\node at (4.19246, 7.31374){\small $\aaa$};
		\node at (5.36168, 9.33185){\small $\bbb$};
		\node at (6.41879, 7.40984){\small $\ccc$};
		\node[inner sep=0.5, draw, shape=circle] at (5.3777, 8.19466){\tiny 1};

		\begin{scope}[xshift=1cm]
				\draw (16.63006, 7.00819)--(16.70929, 3.48237)--(14, 7)--(15.1716, 9.98714)--(17.6239, 9.98714)--(19, 7)--(16.65755, 2.11326)--(16.70929, 3.48237)
				(14, 7)--(19, 7)
				(15.1716, 9.98714)--(16.63006, 7.00819)--(17.6239, 9.98714)
				(14, 7)--(13.46079, 4.43315)--(16.70929, 3.48237);
				
				\node at (15.84606, 9.60214){\small $\aaa$};
				\node at(17.19082, 9.6717){\small $\bbb$};
				\node at (16.61118, 8.00235){\small $\ccc$};
				\node[inner sep=0.5, draw, shape=circle] at (16.51844, 9.09206){\tiny 1};
				
				\node at (16.1011, 7.32997){\small $\aaa$};
				\node at (14.68679, 7.37634){\small $\bbb$};
				\node at (15.24324, 8.99932){\small $\ccc$};
				\node[inner sep=0.5, draw, shape=circle] at (15.24324, 8.09509){\tiny 2};
				
				\node at (18.39646, 7.39953){\small $\aaa$};
				\node at (17.23719, 7.4459){\small $\bbb$};
				\node at (17.70089, 8.99932){\small $\ccc$};
				\node[inner sep=0.5, draw, shape=circle] at (17.74727, 8.04872){\tiny 3};
				
				\node at (16.2634, 6.61123){\small $\aaa$};
				\node at (14.89546, 6.61123){\small $\bbb$};
				\node at (16.2634, 4.5941){\small $\ccc$};
				\node[inner sep=0.5, draw, shape=circle] at(15.68376, 5.91566){\tiny 4};
				
				\node at (18.25734, 6.63441){\small $\aaa$};
				\node at (17.09807, 6.6576){\small $\bbb$};
				\node at (16.98215, 3.62031){\small $\ccc$};
				\node[inner sep=0.5, draw, shape=circle] at(17.3763, 5.45196){\tiny 5};
				
				\node at (13.99123, 4.77958){\small $\aaa$};
				\node at (14.26945, 5.93885){\small $\bbb$};
				\node at (15.40554, 4.57091){\small $\ccc$};
				\node[inner sep=0.5, draw, shape=circle] at(14.43175, 5.2201){\tiny 5};
		\end{scope}
	
	\begin{scope}[xshift=20cm]
	\draw (7, 7)--(7.01336, 3.39949)--(11.03161, 5.13546)--(6.99539, 5.15148)--(3.54613, 7.00203)--(3.55179, 3.40565)--(6.99539, 5.15148)
	(10.02255, 5.13546)--(8.34951, 9.38681)--(7, 7)
	(7,7)--(10,5.15)
	(7.01336, 3.39949)--(3.55179, 3.40565)
	(7, 7)--(3.54613, 7.00203);
	
	\node at (3.87212, 4.03031){\small $\aaa$};
	\node at(3.90416, 6.35273){\small $\bbb$};
	\node at (6.1465, 5.15148){\small $\ccc$};
	\node[inner sep=0.5, draw, shape=circle] at (4.84915, 5.13546){\tiny 1};
	
	\node at (6.75514, 4.65496){\small $\aaa$};
	\node at (6.70709, 3.742){\small $\bbb$};
	\node at (4.86516, 3.75802){\small $\ccc$};
	\node[inner sep=0.5, draw, shape=circle] at (5.82617, 4.03031){\tiny 2};
	
	\node at (6.70709, 6.73714){\small $\aaa$};
	\node at (6.65904, 5.80817){\small $\bbb$};
	\node at (4.76906, 6.75315){\small $\ccc$};
	\node[inner sep=0.5, draw, shape=circle] at (5.92227, 6.35273){\tiny 3};
	
	\node at (7.36377, 3.95022){\small $\aaa$};
	\node at (7.39581, 4.84716){\small $\bbb$};
	\node at (9.63815, 4.86318){\small $\ccc$};
	\node[inner sep=0.5, draw, shape=circle] at (8.19665, 4.60691){\tiny 4};
	
	\node at (7.29971, 5.48783){\small $\aaa$};
	\node at (7.23564, 6.4){\small $\bbb$};
	\node at (9.06155, 5.48783){\small $\ccc$};
	\node[inner sep=0.5, draw, shape=circle] at (7.78021, 6.01638){\tiny 5};
	
	\node at (8.30876, 8.62711){\small $\aaa$};
	\node at (7.55598, 7.18561){\small $\bbb$};
	\node at (9.38188, 5.90427){\small $\ccc$};
	\node[inner sep=0.5, draw, shape=circle] at (8.38885, 7.28171){\tiny 6};
	\end{scope}
		
	\end{tikzpicture}
	\caption{Three possibilities for $\aaa^2\bbb^2\ccc$.}
	\label{naabbc-1}
\end{figure}

\subsubsection*{Case $\aaa^2\bbb\ccc^2$}
From Table \ref{V_5}, this case has a half vertex only if $f=8l(l\ge3)$.

For $f=8l$, we have $\text{AVC}=\{\aaa^2\bbb\ccc^2,\bbb^\frac{f}{4};\,\bbb^\frac{f}{8}\}$. This implies $\aaa$ and $\ccc$ are not in any half vertex. Then we have the AAD $\aaa^2\bbb\ccc^2=\thin^{\ccc}\bbb^{\aaa}\thin^{\bbb}\aaa^{\ccc}\thin^{\ccc}\aaa^{\bbb}\thin^{\bbb}\ccc^{\aaa}\thin^{\aaa}\ccc^{\bbb}\thin,\thin^{\ccc}\bbb^{\aaa}\thin^{\bbb}\aaa^{\ccc}\thin^{\aaa}\ccc^{\bbb}\thin^{\bbb}\aaa^{\ccc }\thin^{\aaa}\ccc^{\bbb}\thin$ or $\thin^{\ccc}\bbb^{\aaa}\thin^{\bbb}\aaa^{\ccc}\thin^{\aaa}\ccc^{\bbb}\thin^{\bbb}\ccc^{\aaa}\thin^{\ccc}\aaa^{\bbb}\thin$.
We have the AAD $\bbb^{\aaa}\thin^{\aaa}\bbb\cdots$ at $\bbb^{\frac{f}{8}}$. This gives a vertex $\aaa^{\bbb}\thin^{\bbb}\aaa\cdots$, a contradiction.

\section{Conclusion}

After the complete classification of non-side-to-side tilings by congruent triangles with any irrational angle, the remaining case of rational triangles (triangles with all angles being rational in degree) will be studied in \cite{jwy}. The techniques developed in this paper, especially the matching lemma and the AAD, remain effective. For example, we explain how the AVC and a new tiling are derived for the triangle with 
\[\aaa=\frac{1}{2},\quad \bbb=\frac{1}{3},\quad \ccc=\frac{1}{6}+\frac{4}{f}, \quad f>24.\]

If $\aaa^m\bbb^n\ccc^l$ is a vertex, then the angle sum of the vertex is (we omit $\pi$ on both sides)
\[\frac{1}{2}m+\frac{1}{3}n+(\frac{1}{6}+\frac{4}{f})l=2\]

This implies
\[(3m+2n+l-12)f=-24l\]

We have $\ccc>\frac{1}{6}$ and the specific values of $\aaa$, $\bbb$. This implies $m\leq4,n\leq6,l\leq12$. We substitute the finitely many combinations of the exponents satisfying these bounds into the equation above and solve for $f$. Those combinations yielding $f\geq30$ are given in Table \ref{rational}. In the table, ``$f=all$" means that the angle combinations can be vertices for any $f$. 

	\begin{table}[H]
	\centering
	\caption{AVC for Case $\aaa=\frac{1}{2},\bbb=\frac{1}{3},\ccc=\frac{1}{6}+\frac{4}{f}$.}
	\begin{tabular}{|c|c|c|c|}
		\hline
		$f$&AVC&
		$f$&AVC\\
		\hline
		all&$\aaa^4$, $\aaa^2\bbb^3$, $\bbb^6$
		&30&$\aaa\ccc^5$\\
		32&$\aaa\bbb\ccc^4$
		&36&$\bbb\ccc^6$, $\aaa\bbb^2\ccc^3$\\
		40&$\bbb^2\ccc^5$
		&48&$\aaa^3\bbb^2$, $\aaa^2\ccc^4$, $\aaa\bbb^3\ccc^2$, $\aaa\ccc^6$, $\bbb^3\ccc^4$, $\ccc^8$\\
		56&$\bbb\ccc^7$
		&60&$\aaa\bbb\ccc^5$\\
		72&$\aaa^2\bbb\ccc^3$, $\bbb^4\ccc^3$, $\bbb^2\ccc^6$, $\ccc^9$
		&84&$\aaa\ccc^7$\\
		96&$\aaa\bbb^2\ccc^4$, $\bbb\ccc^8$
		&120&$\aaa^2\ccc^5$, $\bbb^3\ccc^5$, $\ccc^{10}$\\
		144&$\aaa\bbb\ccc^6$
		&168&$\bbb^2\ccc^7$\\
		192&$\aaa\ccc^8$
		&216&$\bbb\ccc^9$\\
		264&$\ccc^{11}$
		&&\\
		\hline
	\end{tabular}
	
	\label{rational}
\end{table}

\begin{figure}[H]
	\centering
	\includegraphics[width=0.2\textwidth]{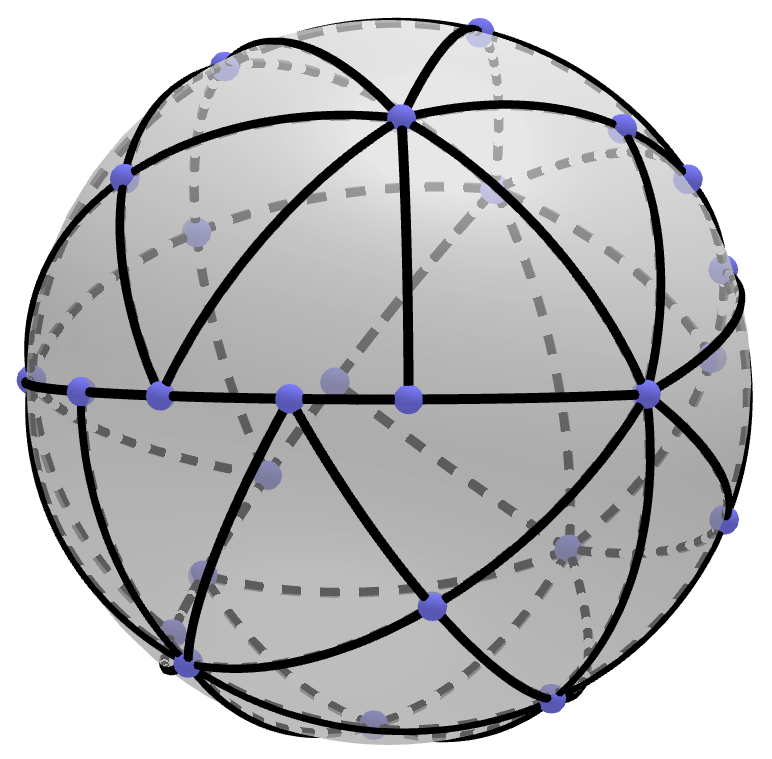}
	\caption{$f=36,T(6\aaa^4,2\bbb^6,6\bbb\ccc^6; \,6\aaa^2, 6\bbb^3)$.}
	\label{yaccccc-3}
\end{figure}

When $f=30$, we find a tiling which are also on Dawson’s website at  \url{https://cs.smu.ca/~dawson}. When $f=36$, we find a new tiling $T(6\aaa^4,2\bbb^6,6\bbb\ccc^6; \,6\aaa^2, 6\bbb^3)$ in Figure \ref{yaccccc-3}. The extended edges are $c+c+a=a
+c+c$. This tiling is neither in \cite{ dawson2006, dawson2006-1, dawson2007}  nor on Dawson’s website.

More details will given in \cite{jwy}.

\end{document}